\documentclass[11pt]{article}

\usepackage{hyperref}
\usepackage{benstyle}
\usepackage{geometry}
\usepackage{enumerate}
\usepackage{longtable}
\usepackage{caption}
\usepackage{comment}
\usepackage{amsmath,pdftexcmds}
\usepackage{enumitem}
\usepackage{microtype}

\newenvironment{foldme}{}{}
\usepackage{cite}

\title{Universal, sample-optimal algorithms for recovery of anisotropic functions from i.i.d. samples}
\author{Ben Adcock \and Avi Gupta\thanks{Department of Mathematics, Simon Fraser University, Canada, \url{avi_gupta@sfu.ca}}\footnotemark[1]}

\begin{document}

\maketitle

\begin{abstract}
A key problem in approximation theory is the recovery of high-dimensional functions from samples. In many cases, the functions of interest exhibit anisotropic smoothness, and, in many practical settings, the nature of this anisotropy may be unknown a priori. Therefore, an important question involves the development of \textit{universal algorithms}, namely, algorithms that simultaneously achieve optimal or near-optimal rates of convergence across a range of different anisotropic smoothness classes. In this work, we consider universal approximation of periodic functions that belong to anisotropic Sobolev spaces and anisotropic dominating mixed smoothness Sobolev spaces. Our first result is the construction of a universal algorithm. This recasts function recovery as a sparse recovery problem for Fourier coefficients and then exploits compressed sensing to yield the desired approximation rates. Note that this algorithm is \textit{nonadaptive}, as it does not seek to learn the anisotropic smoothness of the target function. We then demonstrate optimality of this algorithm up to a dimension-independent polylogarithmic factor. We do this by presenting a lower bound for the adaptive $m$-width for the unit balls of such function classes. Finally, we demonstrate the necessity of nonlinear algorithms. We show that universal linear algorithms can achieve rates that are at best suboptimal by a dimension-dependent polylogarithmic factor. In other words, they suffer from a curse of dimensionality in the rate -- a phenomenon which justifies the necessity of nonlinear algorithms for universal recovery.
\end{abstract}



\section{Introduction}

Many problems require the approximation of functions of $d \gg 1$ that are \textit{anisotropic}, i.e., they exhibit differing degrees of smoothness with respect to different variables. Approximation of high-dimensional, anisotropic functions has been intensively studied over the last several decades. However, the vast majority of existing results consider situations where the anisotropic smoothness is known and strive to design algorithms using this information. This assumption may be unreasonable in practice, particularly in simulation settings where the underlying function is accessible only via a black-box numerical routine. In this paper, we study the significantly more challenging \textit{unknown anisotropy} setting, where algorithms are required to be performant for all possible anisotropic behaviours. We refer to such algorithms as \textit{universal}. Our main contributions are concrete algorithms that achieve near-optimal rates for universal recovery of anisotropic functions. Furthermore, we establish the necessity of nonlinear approximation for universal recovery, by showing that the best rates that can be achieved by linear algorithms are necessarily suboptimal.

\subsection{Problem setting}

We consider $d$-variate functions $f : \bbT^d \rightarrow \bbC$ defined on the $d$-dimensional torus $\bbT^d : = [-\pi,\pi)^d$. Let
\be{
\label{fourier_basis}
\phi_n(x) = \frac{1}{(2\pi)^{d/2}}\E^{\I n \cdot x}, \quad n \in \mathbb{Z}^d
}
denote the $L^2(\bbT^d)$-orthonormal Fourier basis and
\begin{align}\label{fourier_coefficients_notation}
\hat{f}_n = \int_{\bbT^d}f(x) \phi_{-n}(x) \D x,\qquad n \in \bbZ^d,
\end{align}
denote the Fourier coefficients of $f\in L^1(\bbT^d)$. Motivated by the settings studied in \cite{temlyakov2018multivariate,dung2018hyperbolic,moeller2025instance,jahn2023sampling,kuhn2021anisotropic,cobos2016optimal,griebel2014fast} and many others, we consider the following two types of anisotropic Sobolev spaces:

\defn{[Dominating Mixed Smoothness Sobolev Spaces]
\label{def:dom-mixed-space}
  Let $\alpha = (\alpha_1, \ldots, \alpha_d) \in [0,\infty)^d$. We define the \textit{anisotropic dominating mixed smoothness Sobolev space $H^{\alpha}_{\mathsf{mix}}(\bbT^d)$} as 
\begin{align*}
H^{\alpha}_{\mathsf{mix}}(\bbT^d) = \left\{f\in L^2(\bbT^d): \|f\|^2_{H^\alpha_{\mathsf{mix}}} := \sum_{n \in \mathbb{Z}^d} \prod_{j=1}^d (1 + |n_j|)^{2\alpha_j} |\hat{f}_n|^2<\infty\right\}.
\end{align*}
}
\defn{[Anisotropic Sobolev Spaces]
\label{def:anisotropic-space}
Let $\beta = (\beta_1, \ldots, \beta_d) \in [0,\infty)^d$. We define the \textit{anisotropic Sobolev space $H^\beta(\bbT^d)$} as
\begin{align*}
H^\beta(\bbT^d) = \left\{f\in L^2(\bbT^d):\|f\|^2_{H^\beta} := \sum_{n \in \mathbb{Z}^d} \left( 1 + \sum_{j=1}^d |n_j|^{\beta_j} \right)^2 |\hat{f}_n|^2<\infty\right\}.
\end{align*}
}

These two constructions represent two distinct naturally-arising anisotropic behaviors. 
We refer to $\alpha$ and $\beta$ as \textit{anisotropy parameters}, as they control the amount of smoothness in each coordinate. For convenience, we also define
\be{
\label{h-p-def}
h(\alpha)=\min_{\substack{i\in[d]}}\alpha_i,\quad p(\alpha)=|\{j\in[d]:\alpha_j=h(\alpha)\}|,\qquad \alpha \in [0,\infty)^d,
}
and
\be{
\label{g-def}
g(\beta) = \left ( \sum^{d}_{j=1} \frac{1}{\beta_j} \right )^{-1},\qquad \beta \in (0,\infty)^d.
}
As we shall see, these quantities play a key role in the approximation rates.

For either class of spaces, our aim is to design algorithms that achieve optimal approximation rates for \textit{any} value of the anisotropy parameters $\alpha$ or $\beta$. Throughout, our data is in the form of $m$ pointwise samples of an unknown function $f$, i.e.,
\be{
\label{f-samples}
(x_1 , f(x_1)),\ldots,(x_m,f(x_m)).
}
We generally consider \textit{i.i.d. samples}, where the $x_i$ are drawn independently and identically from the uniform probability measure on $\bbT^d$. Therefore, a consequence of this paper is showing that i.i.d.\ samples constitute near-optimal information for recovery in this setting.

\subsection{Contributions}\label{ss:contributions}

We now summarize our main contributions.

\pbk
\textbf{(A) Existence of universal algorithms (Theorems \ref{t:main-prob-mixed}--\ref{t:non-iid-mix},\ref{t:main-prob-sum}--\ref{t:non-iid-sum}).} We show the existence of a reconstruction map $R : (\bbT^d \times \bbC)^m \rightarrow L^2(\bbT^d)$ that takes the samples \ef{f-samples} and produces an approximation whose $L^2$-norm error behaves, up to constants depending on $\alpha$ and $d$, like
\be{
\label{main-rate-mixed}
\left ( \frac{\log^{p(\alpha)-1}(\widetilde{m}) }{\widetilde{m}} \right )^{h(\alpha)},\qquad \text{where } \widetilde{m} = \frac{m}{\log^3(m) \log(\log(m))}
}
for \textit{all} $f \in H^{\alpha}_{\mathsf{mix}}(\bbT^d)$ and \textit{all} $\alpha > 1/2$. For i.i.d.\ samples, this statement holds with high probability for all $\alpha > 1/2$ (Theorem \ref{t:main-prob-mixed}) or in expectation for all $\alpha \in \cA$ (Theorem \ref{t:main-expec-mix}), where $\cA$ is a bounded, but arbitrary subset of $(1/2,\infty)^d$ and $R$ is independent of $\cA$. In Theorem \ref{t:non-iid-mix}, we also show a deterministic result that asserts the existence of sample points $x^{\star}_1,\ldots,x^{\star}_m$ so that this approximation rate holds deterministically for all $f \in H^{\alpha}_{\mathsf{mix}}(\bbT^d)$ and all $\alpha > 1/2$.

In Theorems \ref{t:main-prob-sum}--\ref{t:non-iid-sum} we show identical results for the spaces $H^{\beta}(\bbT^d)$. Specifically, the error behaves, up to constants depending on $\beta$ and $d$, like
\be{
\label{main-rate-sum}
\left( \frac{1}{\widetilde{m}} \right )^{g(\beta)},\qquad \text{where } \widetilde{m} = \frac{m}{\log^3(m) \log(\log(m))}
}
for \textit{all} $f \in H^{\beta}(\bbT^d)$ and \textit{all} $\beta$ for which $g(\beta) > 1/2$.

\pbk
\textbf{(B) Optimality of the rates (Theorems \ref{t:lower-bound-of-adaptive-m-widths-Halpha-mix},\ref{t:lower-bound-of-adaptive-m-widths-Hbeta}).} Our second contribution is to show that the rates \ef{main-rate-mixed} and \ef{main-rate-sum} are optimal up to the polylogarithmic factor appearing in the definition of $\widetilde{m}$. We do this by lower bounding suitable widths for the unit balls of these classes. 

In general, given $1 \leq p \leq \infty$ and $F \subseteq C(\bbT^d)$, we define the \textit{adaptive $m$-width} of $F$ in $L^p$ as
\be{
\label{adaptive-width}
 \varepsilon_m(F,L^p) =  \inf_{\substack{R \in \cR_{\mathsf{nonlin}} \\ S \in \cS_{\mathsf{ada}} }} \sup_{f \in F} \nm{f - R(S(f))}_{L^p},
}
where $\cR_{\mathsf{nonlin}}$ is the set of all maps (linear or nonlinear) $\bbC^m \rightarrow L^p(\bbT^d)$ and $\cS_{\mathsf{ada}}$ is the set of \textit{adaptive} linear maps $C(\bbT^d) \rightarrow \bbC^m$. These maps take the form
\bes{
S(f) = (S_1(f),S_2(f,S_1(f)),\ldots,S_m(f,S_1(f),\ldots,S_{m-1}(f)) ),
}
where $S_1 : C(\bbT^d) \rightarrow \bbC$ is linear and, for $i \geq 2$, $S_i : C(\bbT^d) \times \bbC^{i-1}$ is linear in its first argument.
There are many other types of widths, but these two will suffice for our purposes. 
Specializing to the above spaces, one has the following estimates:
\ea{
\label{lb-contrib}
 \varepsilon_m(B(H^{\alpha}_{\mathsf{mix}}),L^2) \asymp_{\alpha,d} \left ( \frac{\log^{p(\alpha)-1}(m)}{m} \right )^{h(\alpha)} ,
\qquad
 \varepsilon_m(B(H^{\beta},L^2)) \asymp_{\beta,d} \left ( \frac{1}{m} \right )^{g(\beta)}.
}
Note that these estimates are closely related to those found in \cite{byrenheid2017optimal}, although our proof is different -- see \S \ref{s:lower-bounds-widths} for further discussion.
The lower bounds confirm near-optimality of the algorithms constructed in (A). However, these algorithms are fundamentally stronger than those associated with the lower bounds. First, they are universal, meaning that they are independent of $\alpha$ or $\beta$, whereas \ef{lb-contrib} allows the associated algorithms to depend on $\alpha$ or $\beta$. Second, our results use i.i.d.\ sampling from the uniform measure, whereas $\varepsilon_m$ allows arbitrary linear (and adaptive) information, which may not be pointwise samples (standard information) and may also depend on both $m$ and $\alpha$ or $\beta$.

\pbk
\textbf{(C) Necessity of nonlinear algorithms (Theorem \ref{t:main_theorem_index_of_universality}, Corollaries \ref{cor-necessity-Hbeta} and \ref{cor-necessity-Halphamix}).} As noted, the reconstruction maps established in (A) are nonlinear. In our final contribution, we confirm that nonlinearity is necessary. Specifically, in Corollary \ref{cor-necessity-Halphamix} we show the following. Let $\cA$ be any arbitrary subset of $(1/2,\infty)^d$ with nonempty interior and for which $(\alpha_0,\ldots,\alpha_0) \in \mathrm{int}(\cA)$ for some $\alpha_0 > 0$. Now let $G : C(\bbT^d) \rightarrow L^2(\bbT^d)$ be a linear operator of rank at most $m$ for which
\be{
\label{universality-condition}
\sup_{f \in B(H^{\alpha}_{\mathsf{mix}})} \nm{f - G(f)}_{L^2} \lesssim_{\cA,d} n^{-h(\alpha)} (\log(n))^{h(\alpha)(p(\alpha)-1)},\quad \forall \alpha \in \cA,
}
for some $n \in \bbN$. In other words, using $m$ linear (but otherwise arbitrary) samples $G$ achieves the optimal error rate based on $n$ samples over $B(H^{\alpha}_{\mathsf{mix}})$, \textit{uniformly} in $\alpha \in \cA$.
Then, necessarily,
\bes{
n \lesssim_{\cA,d} \frac{m}{(\log m)^{d-1}}.
}
Conversely, picking one of the nonlinear algorithms $H: C(\bbT^d) \rightarrow L^2(\bbT^d)$ from Contribution (A), we see that
\bes{
\sup_{f \in B(H^{\alpha}_{\mathsf{mix}})} \nm{f - H(f)}_{L^2} \lesssim_{\alpha,d} n^{-h(\alpha)} (\log(n))^{h(\alpha)(p(\alpha)-1)},\quad \forall \alpha \in \cA,
}
where $n \in \bbN$ satisfies
\bes{
n \geq \frac{m}{\log^3(m) \log(\log(m))}.
}
Thus, universal linear algorithms incur a curse of dimensionality in the factor $(\log m)^{d-1}$, and nonlinear algorithms are superior to linear algorithms whenever $d > 4$.
In Corollary \ref{cor-necessity-Hbeta} we show the analogous result for the $H^{\beta}$ spaces, with an identical conclusion.

\subsection{Related literature and discussion}

Isotropic and anisotropic spaces have been studied extensively in approximation theory and numerical analysis. Anisotropic Sobolev spaces arise in various applications, including stochastic and parametric PDEs arising in uncertainty quantification, high-dimensional PDEs arising in quantum mechanics, kinetic theory, finance, fluid dynamics and beyond, operator learning and numerous other fields. 
Many works have analyzed various widths corresponding to isotropic and anisotropic Sobolev spaces. See \cite{temlyakov2018multivariate,dung2018hyperbolic,moeller2025instance,krieg2025sampling,dai2026survey} for a modern overview, and, e.g., \cite{kuhn2015approximation,cobos2016optimal,chen2017preasymptotics,kuhn2021anisotropic} for related results on tractability. The more classical results on sampling widths and linear widths in the periodic setting concern the widths of classical isotropic and anisotropic Sobolev spaces. Extending the isotropic results of H\"ollig \cite{hollig1980diameters} to the anisotropic setting, Temlyakov \cite{temlyakov2018multivariate} established sharp asymptotic decay rates for the linear widths of anisotropic Sobolev classes.
Temlyakov \cite{temlyakov2018multivariate} also derives sharp asymptotic decay rates for the linear sampling widths of anisotropic Sobolev classes for all $1 \le p,q \le \infty$, showing that under appropriate smoothness assumptions they decay at the same order as the corresponding linear widths. Here, $q$ characterizes the integrability of the partial derivatives (or, in the case of fractional smoothness, the summability of the Fourier coefficients) and $p$ characterizes the error norm ($L^p$). It is well known (see \cite[Remark~4.21(c)]{novak2008trac}) 
that for classical isotropic Sobolev spaces, nonlinear sampling algorithms do not improve the asymptotic rate compared to linear ones. The asymptotic rates they achieve are listed in \cite[Remark~4.21(c)]{novak2008trac}.

We now turn to anisotropic Sobolev spaces of mixed smoothness. 
For the case $p = q = 2$, \cite{telyakovskii1964some} and \cite{mitjagin1962approximation} established matching asymptotic upper and lower bounds for the $L^2$ linear widths of these classes, yielding sharp asymptotic estimates. Further tractability results in this setting were obtained in \cite{kuhn2021anisotropic}.
More generally, for $1 \le q < p < \infty$ with $p \ge 2$, Byrenheid et al. \cite{byrenheid2017optimal} established sharp asymptotic estimates for the linear widths, sampling widths, and linear sampling widths of anisotropic Sobolev classes of mixed smoothness.
Specific attention has been given to $L^2$ linear sampling widths of functions from Hilbert spaces ($p=q=2$). A breakthrough result by Krieg and M. Ullrich \cite{krieg2021function} showed that sampling recovery for reproducing kernel Hilbert spaces in $L^2$ is asymptotically as powerful as linear approximation. Building on their idea, Dolbeault et al. \cite{dolbeault2023sharp} provided further refinements that can be applied to obtain sharp upper bounds for linear sampling widths of isotropic mixed smoothness Sobolev spaces in $L^2$ (see \cite[section 4.2]{jahn2023sampling}).
Moeller et al. \cite{moeller2025instance} provide upper bounds for nonlinear sampling rates for many regimes of $p$ and $q$, which are shown to be nearly sharp in the case $1<p<2<q<\infty$ with $\frac{1}{p}+\frac{1}{q}>1$. 
Jahn et al. \cite{jahn2023sampling} provide upper bounds for $L^2$ nonlinear sampling widths of Sobolev spaces with mixed smoothness for several regimes of $q$. In the Hilbert space case $p=2$, corresponding to $L^2$ recovery, no improvement in recovery rates can be achieved by nonlinear sampling algorithms over linear ones, in view of \cite[Theorem 4.8]{novak2008trac}. In contrast, for $1<q<2$, they establish that nonlinear sampling methods yield strictly improved decay rates, and they prove that these rates are sharp in this regime. Dai and Temlyakov \cite{dai2024random} investigated sampling recovery for anisotropic Sobolev classes in the $L^2$ norm. In the regime $1<q<2$ (with $p=2$), they derived upper bounds for the sampling widths.
More generally, Kosov and Temlyakov \cite{kosov2025sampling} extended these results to the $L^p$ setting. For $1<q\le 2\le p<\infty$, they established corresponding upper bounds for the sampling widths in $L^p$.

Many other kinds of widths such as Kolmogorov widths,
orthowidths and Gelfand widths have been analyzed (see \cite{temlyakov2018multivariate,dung2018hyperbolic} for an overview). Many of these widths lower bound sampling widths and linear widths and are hence useful tools to show sharpness of the asymptotic upper bounds for linear widths and sampling widths.

Our work differs from this literature in that we consider universal algorithms, delivering near-optimal rates for arbitrary values of the anisotropy parameters. As noted, our algorithms employ compressed sensing. This has been used for over a decade to develop efficient algorithms for high-dimensional function approximation \cite{rauhut2012sparse,rauhut2016interpolation,adcock2022sparse,adcock2024efficient}. More recently, ideas from sparse recovery and the related concepts such as universal discretization \cite{dai2023universal} have been used to derive new results on widths of certain spaces of multivariate functions \cite{adcock2024optimal,adcock2025optimal,jahn2023sampling,moeller2026best,moeller2024high,moeller2025sampling,dai2023universal,dai2025universal,dai2024random,kosov2025sampling}. The majority of results focus on isotropic spaces, such as isotropic mixed smoothness Sobolev spaces and mixed Wiener spaces \cite{moeller2026best,moeller2025instance,moeller2025sampling,krieg2024tractability,kolomoitsev2023sparse,jahn2023sampling,krieg2025sampling,nguyen2022s,moeller2024high,moeller2023gelfand
}. See also \cite{adcock2024optimal,adcock2025optimal,adcock2024optimalb} for results on anisotropic spaces of infinite-dimensional holomorphic functions.

Some works have applied an iterative perspective to approximating functions from high-dimensional periodic Sobolev spaces with unknown anisotropy, where they learn the unknown anisotropy of the function and adapt their approximation procedure accordingly. Bartel and Schröter \cite{bartel2025learning} develop one such approach. Their algorithm alternates between approximating the function and estimating its smoothness from the current approximation and sampling information, with each step feeding into the next to progressively improve both. Our algorithm is, by contrast, \textit{nonadaptive}. It does not strive to estimate the smoothness of the target function, yet it still achieves near-optimal approximations universally without this step.

The notion of universality that we use has also appeared in statistical learning theory, where
it refers to algorithms that achieve optimal rates without a priori knowledge of the
  smoothness of the target function. Binev et al.\cite{binev2005universal,binev2007universal}
  construct such estimators for regression functions from i.i.d.\ samples, achieving
  minimax-optimal rates over the isotropic approximation classes $\mathcal{A}^s$ and
  $\mathcal{B}^s$, which in the Lebesgue measure setting correspond to isotropic Besov
  spaces $B^s_\infty(L_2)$ and $B^s_q(L_\tau)$. In contrast to our noiseless setting, this line of work operates in the noisy
  regression setting, where recovery rates typically take a much slower form than the noiseless case (see, e.g., \cite{devore2025optimal}). Their setting differs from ours in two other ways.
First, their algorithm constructs the
approximation space adaptively from the data, selecting a partition for the domain by thresholding
empirically estimated Haar-like coefficients and approximating by piecewise constant \cite{binev2005universal} or piecewise  polynomial \cite{binev2007universal} functions. Our algorithm is \textit{nonadaptive} by comparison, in that it fixes its structure in advance
and does not iteratively adapt to the data. Second, their results
concern universal recovery in isotropic smoothness classes, whereas we consider anisotropic smoothness classes.

Our focus on unknown anisotropy continues a line of work initiated in \cite{adcock2024optimal,adcock2025optimal,adcock2024optimalb} on infinite-dimensional holomorphic functions that arise frequently in parametric PDEs and operator learning. We pursue a similar approach in this paper -- notably the use of compressed sensing tools -- however, both the setting (Sobolev regularity versus holomorphic regularity) and the techniques employed are quite different. Similar to past works \cite{adcock2022sparse,adcock2024efficient,adcock2025optimal,moeller2026best,moeller2025instance,moeller2025sampling}, our reconstruction map solves a $\ell^1$-minimization based on the so-called \textit{Square-Root LASSO} \cite{belloni2011square-root}. One could also consider greedy algorithms, such as Orthogonal Matching Pursuit (OMP) or weak OMP, as considered in \cite{dai2024random,kosov2025sampling,dai2025universal,moeller2025instance}.

We focus on i.i.d. sampling from the uniform measure -- a situation that is often encountered in practice \cite{adcock2023monte}. Our work therefore contributes to another recent line of research \cite{dai2024random,sonnleitner2023power,adcock2024optimal,adcock2025optimal,krieg2022recovery,krieg2024random,krieg2025function,dai2024random} that determines scenarios where i.i.d. sampling is near optimal. As demonstrated by the lower bounds, our results are optimal up to the polylogarithmic term $\log^3(m) \log(\log(m))$ appearing in \ef{main-rate-mixed}-\ef{main-rate-sum}. See Remark \ref{rem:log-term} for further discussion on this factor.

Contribution (C) of our paper is based on Temlyakov's work in \cite[\S 5.4]{temlyakov2018multivariate} and \cite{temlyakov1988approximation}, which considered the spaces $H^{\beta}$ only. We extend these results to the $H^{\alpha}_{\mathsf{mix}}$ spaces. Note that \cite{jahn2023sampling} compares linear and nonlinear widths for a variety of Wiener type spaces, as well as mixed smoothness Sobolev spaces, and also shows the superiority of nonlinear algorithms in various settings.
Jahn et al.~\cite{jahn2023sampling} show that for mixed Wiener spaces  
nonlinear recovery strictly outperforms linear recovery in $L^2$ (see \cite[Remark~4.5(i)]{jahn2023sampling}). 
In contrast, for $L^2$ recovery in mixed smoothness Sobolev spaces
 (the Hilbert case $q=2$) no such improvement occurs (see \cite[Theorem~4.8]{novak2008trac}), whereas for general mixed smoothness Sobolev spaces, nonlinear recovery decays faster than linear widths for $L^2$ recovery (see \cite[Remark~4.17]{jahn2023sampling}). Moeller et al. \cite{moeller2025instance} extend the conclusion in \cite[Remark~4.5(i)]{jahn2023sampling}, that we stated above, to the case of 
mixed Wiener spaces with $2\le p<\infty$; that is, they show that
nonlinear sampling achieves a strictly faster main rate than linear methods in this setting (see \cite[Remark~4.6]{moeller2025instance}). However, they do not consider universality across ranges of parameters and the necessity of nonlinear algorithms for universal recovery.

Many of the aforementioned works consider approximation in $L^p$ spaces for general $p \in [1,\infty]$. For clarity of exposition, we have chosen to consider $p = 2$ only. However, we anticipate that many of our main results could be extended. This is a question for future work.

\subsection{Outline}

The outline of the remainder of this paper is as follows. We commence in \S \ref{s:best-s-term} by deriving best $s$-term approximation rates in the Fourier basis for the $H^{\alpha}_{\mathsf{mix}}$ and $H^{\beta}$, estimates that will be crucial later. In \S \ref{s:universal-algorithms-thms}-\ref{s:universal-lower} we address contributions (A), (B) and (C) of \S \ref{ss:contributions}, respectively.

\section{Best $s$-term approximation in $H^{\alpha}_{\mathsf{mix}}$ and $H^{\beta}$}\label{s:best-s-term}

The main results in this paper exploit the concept of best $s$-term approximation in the Fourier basis. Therefore, in this section, we establish a series of estimates for the best $s$-term approximation of functions in $H^{\alpha}_{\mathsf{mix}}$ and $H^{\beta}$. 
Let $f \in L^2(\bbT^d)$, and let $\{ \phi_n \}_{n \in \bbZ^d}$ and $\hat{f}_n$ ($n \in \bbZ^d$) be as in \ef{fourier_basis} and \ef{fourier_coefficients_notation}, respectively. We write the Fourier series of $f$ as
\begin{align*}
f=\sum_{n \in \mathbb{Z}^d}\hat{f}_n\phi_n
\end{align*}
(with convergence in $L^2$) and recall that Parseval's identity holds, i.e., $\int_{\bbT^d}|f(x)|^2 \D x=\sum_{n\in\mathbb{Z}^d}|\hat{f}_n|^2$.

\subsection{Sequence spaces and best $s$-term approximation}

We first require a series of definitions.

\defn{[$\ell^p$ space]
Let $\Lambda$ be countable and $p > 0$. The \textit{$\ell^p(\Lambda)$ space} is the set
\begin{align*}
    \ell^p(\Lambda) = \left\{ c = (c_\lambda)_{\lambda \in \Lambda}  : 
    \|c\|_p := 
    \begin{cases}
        \left( \sum_{\lambda \in \Lambda} |c_\lambda|^p \right)^{1/p}, & 0 < p < \infty, \\
        \sup_{\lambda \in \Lambda} |c_\lambda|, & p = \infty
    \end{cases}
     < \infty \right\}.
\end{align*}
}

\defn{[Weak $\ell^p$ space]
Let $N := |\Lambda| \in \mathbb{N} \cup \{\infty\}$ and $(c_i^*)_{i=1}^N$ denote a non-increasing rearrangement (by absolute value) of $c=(c_\lambda)_{\lambda \in \Lambda}$, i.e., $|c_i^*| \ge |c_{i+1}^*|$ for all $i = 1,\ldots,N$. For $p > 0$, the \textit{weak $\ell^p(\Lambda)$ space (w$\ell^p$)} is the set
\begin{align*}
    \mathrm{w}\ell^p(\Lambda) = 
    \left\{ c \in \mathbb{C}^\Lambda : 
    \|c\|_{p,\infty} := 
    \begin{cases}
        \sup_{1 \le i \le N} |c_i^*|  i^{1/p}, & 0 < p < \infty, \\
        \sup_{1 \le i \le N} |c_i^*|, & p = \infty
    \end{cases}
    < \infty \right\}.
\end{align*}
}
\defn{[Weak Lorentz $\ell^{p,a}$ space]
Let $N$ and $(c_i^*)_{i=1}^N$ be as in the previous definition. For constants $p>0$ and $a\geq0$, the \textit{weak Lorentz $\ell^{p,a}(\Lambda)$ space (w$\ell^{p,a}$)} is
\begin{align*}
    \mathrm{w}\ell^{p,a}(\Lambda) = 
    \left\{ c \in \mathbb{C}^\Lambda :
    \|c\|_{p,a,\infty} := 
    \begin{cases}
        \sup_{1 \le i \le N} |c_i^*| \, i^{1/p} \bigl(\max\{1,\log i\}\bigr)^{-a/p}, 
        & 0 < p < \infty, \\
        \sup_{1 \le i \le N} |c_i^*|, 
        & p = \infty
    \end{cases}
    < \infty \right\}.
\end{align*}
}

Note that these spaces satisfy the relations $\ell^p \subsetneq w \ell^p = w \ell^{p,0} \subsetneq w \ell^{p,a}$ for all $a > 0$.

\defn{[Best $s$-term approximation error]
\label{def:sigma_s}
Let $q > 0$, $\Lambda$ be countable and $z\in \ell^q(\Lambda)$. 
The \emph{$\ell^q$-norm best $s$-term approximation error} of $z$ is defined as
\[
\sigma_s(z)_q = \inf_{\substack{w \in \ell^q(\Lambda) \\ |\mathrm{supp}(w)| \le s}} \|z - w\|_q.
\]

}

Here and elsewhere, $\mathrm{supp}(w) = \{ \lambda \in \Lambda : w_\lambda \neq 0 \}$. If $z^* = (z_i^*)^{N}_{i=1}$ denotes a non-increasing rearrangement of $z$, then we can write this as
\[
\sigma_s(z)_q = \Bigg( \sum_{i > s} |z_i^*|^q \Bigg)^{1/q}.
\]

\subsection{Stechkin's inequalities}

The following result is standard in the case $a=0$ and $q<\infty$. See, e.g., \cite[Lem. SM2.1]{adcock2023monte}.
Our lemma below is inspired by this reference. It extends the case $a=0$ to general $a\ge0$, and we also include the case $q=\infty$.
We also allow $s=0$. It should be noted that 
results such as the one below
were never published in this form by Stechkin. See \cite[§7.4]{dung2018hyperbolic} for further details. Nevertheless, we continue to refer to it as \textit{Stechkin's inequality}, as this terminology is widely used in the high-dimensional approximation and compressed sensing communities.

\lem{[Stechkin's inequality in $w\ell^{p,a}$]\label{l:stech2}
Let $0<p<q\le\infty$, $a\ge0$ and $c=(c_i)_{i\in\Lambda}\in\bbC^{|\Lambda|}$, where $\Lambda\subset\bbN$ is an index set and $N=|\Lambda|\in\bbN\cup\{\infty\}$.
Then $c\in w\ell^{p,a}(\Lambda)$ if and only if there exists a constant $C>0$ depending only on $p,q,a,$ and $c$ such that
\begin{align*}
\sigma_s(c)_q \le C \bigl(\max\{1,s\}\bigr)^{\frac1q-\frac1p}
\bigl(\max\{1,\log(\max\{1,s\})\}\bigr)^{\frac{a}{p}}
\end{align*}
for all $s\in\bbN\cup\{0\}$ with $s<N$ when $N<\infty$.
Specifically, if $c\in w\ell^{p,a}(\Lambda)$ and $q<\infty$, then
\begin{align*}
\sigma_s(c)_q \lesssim_{p,q,a}
\|c\|_{p,a,\infty}
\bigl(\max\{1,s\}\bigr)^{\frac1q-\frac1p}
\bigl(\max\{1,\log(\max\{1,s\})\}\bigr)^{\frac{a}{p}}
\end{align*}
for all $s\in\bbN\cup\{0\}$ with $s<N$ when $N<\infty$.
If $c\in w\ell^{p,a}(\Lambda)$ and $q=\infty$, then
\begin{align*}
\sigma_s(c)_\infty \lesssim_{p,a}
\|c\|_{p,a,\infty}
\bigl(\max\{1,s\}\bigr)^{-1/p}
\bigl(\max\{1,\log(\max\{1,s\})\}\bigr)^{a/p}
\end{align*}
for all $s\in\bbN\cup\{0\}$ with $s<N$ when $N<\infty$.
Conversely, if $q<\infty$ and
\begin{align*}
\sigma_s(c)_q \le C \bigl(\max\{1,s\}\bigr)^{\frac1q-\frac1p}
\bigl(\max\{1,\log(\max\{1,s\})\}\bigr)^{\frac{a}{p}}
\end{align*}
for all $s\in\bbN\cup\{0\}$ with $s<N$ when $N<\infty$, then
$
\|c\|_{p,a,\infty} \le 2^{\frac1p+\frac1q}C.
$
If $q=\infty$ and
\begin{align*}
\sigma_s(c)_\infty \le C \bigl(\max\{1,s\}\bigr)^{-1/p}
\bigl(\max\{1,\log(\max\{1,s\})\}\bigr)^{a/p}
\end{align*}
for all $s\in\bbN\cup\{0\}$ with $s<N$ when $N<\infty$, then
$
\|c\|_{p,a,\infty} \le 2^{1/p}C.
$
}

\prf{
Suppose first that $c \in w \ell^{p,a}(\Lambda)$ and let $(c_i^*)_{i=1}^N$ be a nonincreasing rearrangement of $c$.
For $i\in\bbN$ with $i\le N$ when $N<\infty$, we have
\begin{align*}
|c_i^*|\le \|c\|_{p,a,\infty}\, i^{-1/p}\bigl(\max\{1,\log i\}\bigr)^{a/p}.
\end{align*}
We first consider the case $q<\infty$.
Let $s\in\bbN\cup\{0\}$ with $s<N$ when $N<\infty$.
We set $s_0=\max\{3,\left\lceil \E^{\frac{2aq}{q-p}}\right\rceil\}$.
Suppose first that $s\ge s_0$. Then $s\ge3$ and
\begin{align*}
(\sigma_s(c)_q)^q
=\sum_{i>s}|c_i^*|^q
\le \|c\|_{p,a,\infty}^q\sum_{i>s} i^{-q/p}\bigl(\max\{1,\log i\}\bigr)^{aq/p}
=\|c\|_{p,a,\infty}^q\sum_{i>s} i^{-q/p}(\log i)^{aq/p}.
\end{align*}
Since $s\ge3$ and $s\ge \E^a$, the function $x\mapsto x^{-q/p}(\log x)^{aq/p}$ is decreasing for $x\ge s$. Hence
\begin{align*}
(\sigma_s(c)_q)^q\le \|c\|_{p,a,\infty}^q\int_s^\infty x^{-q/p}(\log x)^{aq/p}\D x.
\end{align*}
With the changes of variables $x=\E^t$ and $u=\frac{q-p}{p}t$, the integral becomes
\begin{align*}
\int_s^\infty x^{-q/p}(\log x)^{aq/p}\D x=\int_{\log s}^\infty \E^{(p-q)t/p} t^{aq/p}\D t =\left(\frac{p}{q-p}\right)^{\frac{aq}{p}+1}\int_{\frac{q-p}{p}\log s}^\infty u^{\frac{aq}{p}}\E^{-u}\D u.
\end{align*}
Recall that the \textit{upper incomplete Gamma function} is defined as $\Gamma(r,x)=\int_x^\infty u^{r-1}\E^{-u}\D u$ for $r>0$, $x\ge0$.
We deduce that
\begin{align*}
(\sigma_s(c)_q)^q \le \|c\|_{p,a,\infty}^q\left(\frac{p}{q-p}\right)^{\frac{aq}{p}+1}\Gamma\left(\frac{aq}{p}+1,\frac{q-p}{p}\log s\right).
\end{align*}
We shall use the bound $\Gamma(r,x)\le 2x^{r-1}\E^{-x}$ for $r\ge1$ and $x\ge 2(r-1)$. To see why this holds, let $r\ge1$ and $x>0$. Integration by parts gives
\begin{align*}
\Gamma(r,x)  =x^{r-1}\E^{-x}+(r-1)\int_x^\infty u^{r-2}\E^{-u}\D u
&\le x^{r-1}\E^{-x}+\frac{r-1}{x}\int_x^\infty u^{r-1}\E^{-u}\D u
\\
&=x^{r-1}\E^{-x}+\frac{r-1}{x}\Gamma(r,x).
\end{align*}
If $x\ge 2(r-1)$, then $1-\frac{r-1}{x}\ge \frac12$, which implies that $\Gamma(r,x)\le 2x^{r-1}\E^{-x}$, as claimed.
Applying this with $r=\frac{aq}{p}+1$ and $x=\frac{q-p}{p}\log s$, and using the fact that $s\ge s_0$, we get
\begin{align*}
(\sigma_s(c)_q)^q
&\le \|c\|_{p,a,\infty}^q\left(\frac{p}{q-p}\right)^{\frac{aq}{p}+1}2\left(\frac{q-p}{p}\log s\right)^{\frac{aq}{p}}\E^{-\frac{q-p}{p}\log s}
\\
&=\|c\|_{p,a,\infty}^q 2\left(\frac{p}{q-p}\right)(\log s)^{\frac{aq}{p}} s^{1-\frac{q}{p}}.
\end{align*}
Taking the $q$th root gives
\begin{align*}
\sigma_s(c)_q \lesssim_{p,q,a} \|c\|_{p,a,\infty} s^{\frac1q-\frac1p}(\log s)^{\frac{a}{p}},\quad \forall s \geq s_0.
\end{align*}
We now consider the case $s<s_0$.We have
\begin{align*}
\sigma_s(c)^q_q \leq \sigma_{0}(c)^q_q = \|c\|_q^q=\sum_{i\ge1}|c_i^*|^q=\sum_{1\le i\le s_0}|c_i^*|^q+\sum_{i>s_0}|c_i^*|^q.
\end{align*}
For the finite sum, we use the fact that $|c_i^*|\le \|c\|_{p,a,\infty} i^{-1/p}\bigl(\max\{1,\log i\}\bigr)^{a/p}$ to get
\begin{align*}
\sum_{1\le i\le s_0}|c_i^*|^q
\le \|c\|_{p,a,\infty}^q\sum_{1\le i\le s_0} i^{-q/p}\bigl(\max\{1,\log i\}\bigr)^{aq/p}
\lesssim_{p,q,a}\|c\|_{p,a,\infty}^q.
\end{align*}
Here, we also used the fact that $s_0$ depends on $p$, $q$ and $a$ only.
For the tail, we have $\sum_{i>s_0}|c_i^*|^q=(\sigma_{s_0}(c)_q)^q$.
The bound proved above for $s\ge s_0$ applies when $s=s_0$. This gives $(\sigma_{s_0}(c)_q)^q \lesssim_{p,q,a}\|c\|_{p,a,\infty}^q$.
Combining these estimates yields $\|c\|_q\lesssim_{p,q,a}\|c\|_{p,a,\infty}$. Hence
\begin{align*}
\sigma_s(c)_q\le \|c\|_q\lesssim_{p,q,a}\|c\|_{p,a,\infty},\quad \forall s < s_0.
\end{align*}
Moreover, for $s<s_0$ we have $1\le \max\{1,s\}\le s_0$ and $1\le \max\{1,\log(\max\{1,s\})\}\le \max\{1,\log s_0\}$.
Since $\frac1q-\frac1p<0$, it follows that $(\max\{1,s\})^{\frac1q-\frac1p}\ge s_0^{\frac1q-\frac1p}$,
and therefore
\begin{align*}
s_0^{\frac1p-\frac1q}\bigl(\max\{1,s\}\bigr)^{\frac1q-\frac1p}\bigl(\max\{1,\log(\max\{1,s\})\}\bigr)^{a/p}
\ge 1.
\end{align*}
And since $s_0$ only depends on $p,q$ and $a$, after adjusting constants, we may write
\begin{align*}
\sigma_s(c)_q\lesssim_{p,q,a} \|c\|_{p,a,\infty} \bigl(\max\{1,s\}\bigr)^{\frac1q-\frac1p}\bigl(\max\{1,\log(\max\{1,s\})\}\bigr)^{a/p},\quad \forall s < s_0.
\end{align*}
Combining this with the case $s\ge s_0$ gives
\begin{align*}
\sigma_s(c)_q\lesssim_{p,q,a}\|c\|_{p,a,\infty}\bigl(\max\{1,s\}\bigr)^{\frac1q-\frac1p}\bigl(\max\{1,\log(\max\{1,s\})\}\bigr)^{\frac{a}{p}},\quad \forall s \in \bbN \cup \{0 \},
\end{align*}
as required. 
We now consider the case $q=\infty$. 
Let $s\in\bbN\cup\{0\}$ with $s<N$ when $N<\infty$ and set $t=\max\{1,s\}$. 
Then 
\begin{align*} 
\sigma_s(c)_\infty=\sup_{i>s}|c_i^*|\ = |c_{s+1}^*|\le \|c\|_{p,a,\infty}(s+1)^{-1/p}\bigl(\max\{1,\log(s+1)\}\bigr)^{a/p}. 
\end{align*} 
Since $(s+1)^{-1/p}\le t^{-1/p}$ and $\max\{1,\log(s+1)\}\le 2\max\{1,\log t\}$, we get
\begin{align*} 
\sigma_s(c)_\infty\le 2^{a/p}\|c\|_{p,a,\infty}t^{-1/p}\bigl(\max\{1,\log t\}\bigr)^{a/p},
\end{align*}
as required.

We now prove the converse.
Assume first that $q<\infty$ and that there exists $C>0$ such that for all $s\in\bbN\cup\{0\}$ with $s<N$ when $N<\infty$, 
\begin{align*} 
\sigma_s(c)_q \le C \bigl(\max\{1,s\}\bigr)^{\frac1q-\frac1p}\bigl(\max\{1,\log(\max\{1,s\})\}\bigr)^{\frac{a}{p}}. 
\end{align*} 
Taking $s=0$ gives $\sigma_0(c)_q=\|c\|_q\le C$. Hence $|c_1^*|\le C$. 
Let $n\in\bbN$ with $2n\le N$ when $N<\infty$.
Since $(c_i^*)$ is nonincreasing, we have $|c_i^*|\ge |c_{2n}^*|$ for $n+1\le i\le 2n$. Therefore
\begin{align*} 
n|c_{2n}^*|^q\le \sum_{i=n+1}^{2n}|c_i^*|^q\le \sum_{i>n}|c_i^*|^q=(\sigma_n(c)_q)^q.
\end{align*} 
Using the assumed bound for $s=n$, we get
\begin{align*} 
|c_{2n}^*|\le C\, n^{-1/p}\bigl(\max\{1,\log n\}\bigr)^{a/p}.
\end{align*} 
Now let $n\in\bbN$ with $2n+1\le N$ when $N<\infty$. 
Since $|c_i^*|\ge |c_{2n+1}^*|$ for $n+2\le i\le 2n+1$, we have 
\begin{align*} 
n|c_{2n+1}^*|^q\le \sum_{i=n+2}^{2n+1}|c_i^*|^q\le \sum_{i>n+1}|c_i^*|^q=(\sigma_{n+1}(c)_q)^q. 
\end{align*} 
Using the assumed bound for $s=n+1$ and the fact that $(n+1)/n\le 2$, we obtain 
\begin{align*} 
|c_{2n+1}^*|\le 2^{1/q}C\, (n+1)^{-1/p}\bigl(\max\{1,\log(n+1)\}\bigr)^{a/p}. 
\end{align*} 
Combining these bounds and using the fact that $(2n)^{1/p}\le 2^{1/p}n^{1/p}$ and $(2n+1)^{1/p}\le 2^{1/p}(n+1)^{1/p}$ gives 
\begin{align*} 
\|c\|_{p,a,\infty} & 
\le \sup_{n\ge 1}
\big\{
(2n)^{1/p}\bigl(\max\{1,\log (2n)\}\bigr)^{-a/p}|c_{2n}^*|,
(2n+1)^{1/p}\bigl(\max\{1,\log (2n+1)\}\bigr)^{-a/p}|c_{2n+1}^*|
\big\}
\\
&\le2^{\frac1p+\frac1q}C. 
\end{align*} 
Finally, assume that $q=\infty$ and that there exists $C>0$ such that for all $s\in\bbN\cup\{0\}$ with $s<N$ when $N<\infty$, 
\begin{align*} 
\sigma_s(c)_\infty \le C \bigl(\max\{1,s\}\bigr)^{-1/p}\bigl(\max\{1,\log(\max\{1,s\})\}\bigr)^{a/p}. 
\end{align*} 
Taking $s=0$ gives $|c_1^*|\le C$. 
Let $i\in\bbN$ with $2\le i\le N$ when $N<\infty$ and take $s=i-1$. 
Then $\sigma_{i-1}(c)_\infty=\sup_{j\ge i}|c_j^*|=|c_i^*|$, so 
\begin{align*} 
|c_i^*|\le C (i-1)^{-1/p}\bigl(\max\{1,\log(i-1)\}\bigr)^{a/p}. 
\end{align*} 
Since $(i-1)^{-1/p}\le 2^{1/p} i^{-1/p}$ and $\max\{1,\log(i-1)\}\le \max\{1,\log i\}$ for $i\ge 2$, we get 
\begin{align*} 
|c_i^*|\le 2^{1/p} C  i^{-1/p}\bigl(\max\{1,\log i\}\bigr)^{a/p}. 
\end{align*} 
Therefore $\|c\|_{p,a,\infty}\le 2^{1/p}C $, as required.
} 

\subsection{Best $s$-term approximation rates in $H^{\alpha}_{\mathsf{mix}}$}

We now introduce some notation that will be used in the subsequent results. Given $r > 1$ and $\alpha=(\alpha_1,\alpha_2,\dots,\alpha_d)\in(0,\infty)^d$, let 
\begin{align*}
A(d,r,\alpha) = \left|\left\{ n\in\mathbb{Z}^d : \prod_{j=1}^d (1 + |n_j|)^{\alpha_j} < r\right\}\right|,\qquad a(d,r,\alpha) = \left|\left\{ n\in\mathbb{Z}_+^d : \prod_{j=1}^d (1 + n_j)^{\alpha_j} < r\right\}\right|.
\end{align*}
Observe that 
\begin{align}\label{A_equiv_a}
a(d,r,\alpha)\leq A(d,r,\alpha)\leq2^d\times a(d,r,\alpha).
\end{align}
This lemma, as well as Lemma \ref{lem:cdm} later, is well-established in the literature. See, for example, \cite{chen2017preasymptotics,kuhn2021anisotropic}. We offer short proofs for completeness.
\lem{[]\label{l:adm}
Let $\alpha=(\alpha_1,\alpha_2,\dots,\alpha_d)\in(0,\infty)^d$ and $h(\alpha)$ and $p(\alpha)$ be as in \ef{h-p-def}. 
Then
\begin{align}\label{aa}
a(d,r,\alpha)\asymp_{d,\alpha} r^\frac{1}{h(\alpha)}(\log r)^{p(\alpha)-1},
\end{align}
and
\begin{align}\label{AA}
A(d,r,\alpha)\asymp_{d,\alpha} r^\frac{1}{h(\alpha)}(\log r)^{p(\alpha)-1}.
\end{align}
}
\prf{
Notice that \ef{AA} follows immediately from \ef{aa} and \ef{A_equiv_a}. Hence, we only need to show \ef{aa}. 
To do this,
we fix $p(\alpha)=p\in\bbN$ for arbitrary $\alpha \in (0,\infty)^d$, where $d \in \mathbb{N}$ is unspecified and proceed by induction on $d \geq p(\alpha)$. 
First, for the base case $d=p(\alpha)$, we use the fact that $a(d,r,\alpha)=a(d,r^\frac{1}{h(\alpha)},1_d)$, where $1_d$ is the vector of ones of length $d$. \cite[Lemma B.3]{adcock2022sparse} gives that
\begin{align*}
\frac{1}{(d-1)!}\frac{(\lceil r \rceil-1)(\log (\lceil r\rceil-1))^{d}}{\log (\lceil r\rceil-1) + d}< a(d,r,1_d)<(\lceil r \rceil -1)\frac{(\log(\lceil r \rceil-1 )+d\log2)^{d-1}}{(d-1)!}
\end{align*}
where the lower bound holds for $r>1$ and the upper bound holds for $r \geq r^*$, where $r^*$ depends on $d$ only.
Thus $a(d,r,1_d)\asymp_d r(\log r)^{d-1}$, 
and hence
\begin{align*}
a(d,r,\alpha)=a(d,r^\frac{1}{h(\alpha)},1_d)\asymp_d r^\frac{1}{h(\alpha)}(\log r)^{p(\alpha)-1},\quad d = p(\alpha),
\end{align*}
as required. Now assume that $a(d-1,r,\alpha)\asymp_{d,\alpha}\frac{1}{(p-1)!} r^\frac{1}{h(\alpha)}(\log r)^{p-1}$ for some $d > p(\alpha)$. Without loss of generality, we may assume that $\alpha$ is nondecreasing. Now observe that
\begin{align*}
a(d,r,\alpha) = \sum_{j=1}^{\lfloor r^{1/\alpha_d}\rfloor} a(d-1, r/j^{\alpha_d},\alpha'),\quad \text{where }\alpha' = (\alpha_1,\ldots,\alpha_{d-1}).
\end{align*}
This implies that
\begin{align*}
a(d,r,\alpha) \asymp_{d,\alpha} \sum_{j=1}^{\lfloor r^{1/\alpha_d}\rfloor} (r/j^{\alpha_d})^{1/h(\alpha)} (\log (r/j^{\alpha_d}))^{p-1} .
\end{align*}
We now establish upper and lower bounds separately. For the former, observe that, by assumption, $p(\alpha) < d$, and hence $\alpha_{d} > h(\alpha)$. Therefore
\begin{align*}
a(d,r,\alpha)
&\lesssim_{d,\alpha}  r^{1/h(\alpha)} (\log r)^{p(\alpha)-1}
\sum_{j=1}^{\lfloor r^{1/\alpha_d}\rfloor} j^{-\alpha_{d} / h(\alpha)}
\lesssim_{d,\alpha} r^{1/h(\alpha)} (\log r)^{p(\alpha)-1},
\end{align*}
since the series $\sum_{j=1}^{\infty} j^{-\alpha_{d} / h(\alpha)}$ converges for $\alpha_{d} > h(\alpha)$.
Hence the upper bound in \ef{aa} holds for the given $d > p(\alpha)$.

For the lower bound, we consider the $j = 1$ term only, which gives
\begin{align*}
a(d,r,\alpha)
&\gtrsim_{d,\alpha} (r/1^{\alpha_{d}})^{1/h(\alpha)} (\log (r/1^{\alpha_{d}}))^{p(\alpha)-1}
\gtrsim_{d,\alpha} r^{1/h(\alpha)} (\log r)^{p(\alpha)-1},
\end{align*}
as required.
}

We are now ready to present a result on best $s$-term approximation in $H^{\alpha}_{\mathsf{mix}}$.

\thm{[Best s-term approximation rate in $H^{\alpha}_{\mathsf{mix}}$]\label{t:decay_prod}
Let $\alpha=(\alpha_1,\alpha_2,\dots,\alpha_d)$ and $h(\alpha)$ and $p(\alpha)$ be as in \ef{h-p-def}.
Let $f \in H^{\alpha}_{\mathsf{mix}}$ and write $c = (\hat{f}_n)_{n \in \bbZ^d}$, where $\hat{f}_n$ is as in \ef{fourier_coefficients_notation}. Then
\begin{align*}
\sigma_s(c)_q\lesssim_{\alpha,q} s^{\frac{1}{q}-h(\alpha)-\frac{1}{2}}(\log s)^{h(\alpha)(p(\alpha)-1)}\|f\|_{H^{\alpha}_{\mathsf{mix}}},\quad \forall s \in \bbN,
\end{align*}
for all $q>\frac{1}{h(\alpha)+1/2}$. Moreover, for $q=\infty$,
\begin{align*}
\sigma_s(c)_\infty\lesssim_{\alpha} s^{-h(\alpha)-\frac{1}{2}}(\log s)^{h(\alpha)(p(\alpha)-1)}\|f\|_{H^{\alpha}_{\mathsf{mix}}},\quad \forall s \in \bbN.
\end{align*}
}
\prf{
We first prove the result for $q=2$. 
Given $s\in\bbN$, let $\varepsilon > 0$ be such that $A(d,1/\varepsilon,\alpha) \leq s$, where $A(d,1/\varepsilon,\alpha)$ is as in \ef{AA} with $r = 1/\varepsilon$. Then
\be{
\label{eq:tail_bound_product_sobolev}
\sigma_s(c)^2_2 \leq  \sum_{n : \prod^{d}_{j=1} (1+|n_j|)^{\alpha_j} \geq 1/\varepsilon} | \hat{f}_n |^2
 \leq \varepsilon^2 \sum_{n \in \bbZ^d} \prod^{d}_{j=1} (1+|n_j|)^{2\alpha_j} | \hat{f}_n |^2
 = \varepsilon^2 \nm{f}^2_{H^{\alpha}_{\mathsf{mix}} }.
}
By Lemma \ref{l:adm}, there are constants $c_1 = c_1(d,\alpha),c_2 = c_2(d,\alpha)>0$ with $c_2>(p(\alpha)-1)!$ such that
\begin{align*}
\frac{c_1}{(p(\alpha)-1)!}r^{1/h(\alpha)}(\log(r))^{p(\alpha)-1}\leq A(d,r,\alpha)\leq\frac{c_2}{(p(\alpha)-1)!}r^{1/h(\alpha)}(\log(r))^{p(\alpha)-1}.
\end{align*}
There exists, for a large enough $s$, a positive real number $\varepsilon<1/\E$ such that 
\begin{align}\label{ubHA}
s\leq\frac{c_2}{(p(\alpha)-1)!}\varepsilon^{-1/h(\alpha)}(\log(1/\varepsilon))^{p(\alpha)-1}<s+1,
\end{align}
which immediately implies that $A(d,1/\varepsilon,\alpha) \leq s$, and therefore \ef{eq:tail_bound_product_sobolev} holds. We now upper bound $\varepsilon$.
Taking natural logarithms on both sides of \ef{ubHA} gives
\begin{align*}
\log \left(\frac{c_2}{(p(\alpha)-1)!}\right) 
+ \frac{1}{h(\alpha)} \log \left(\frac{1}{ \varepsilon}\right) 
+ (p(\alpha)-1) \log \log \left(\frac{1}{\varepsilon}\right) 
< \log(s+1).
\end{align*}
Since $\log \left(\frac{c_2}{(p(\alpha)-1)!}\right) \ge 0$ and $(p(\alpha)-1) \log \log \left(\frac{1}{\varepsilon}\right) \ge 0$, we get $\frac{1}{h(\alpha)} \log \left(\frac{1}{ \varepsilon}\right) < \log(s+1)$.
Substituting this bound into \ef{ubHA} leads to
\begin{align*}
s\leq \frac{c_2}{(p(\alpha)-1)!} \varepsilon^{-1/h(\alpha)} h(\alpha)^{p(\alpha)-1} \big(\log(s+1)\big)^{p(\alpha)-1}.
\end{align*}
Rearranging,
we obtain
\begin{align}
\varepsilon\leq \left(h(\alpha)^{p(\alpha)-1}\frac{c_2}{(p(\alpha)-1)!} \right)^{h(\alpha)}s^{-h(\alpha)}\big(\log(s+1)\big)^{h(\alpha)(p(\alpha)-1)}.\label{epsiloninequality}
\end{align}
Combining this with \ef{eq:tail_bound_product_sobolev}, we deduce that
\begin{align*}
\sigma_s(c)_2\lesssim_{\alpha} s^{-h(\alpha)}(\log s)^{h(\alpha)(p(\alpha)-1)}\|f\|_{H^{\alpha}_{\mathsf{mix}}},
\end{align*}
which gives the result for $q = 2$.

Next we consider $q \geq 2$. The previous result and \cf{l:stech2} imply that $c\in w\ell^{\frac{1}{h(\alpha)+1/2},\frac{h(\alpha)(p(\alpha)-1)}{h(\alpha)+1/2}}$.
Reapplying \cf{l:stech2} for $q>\frac{1}{h(\alpha)+1/2}$ now gives
\begin{align*}
\sigma_s(c)_q\lesssim_{\alpha,q} s^{\frac1q-h(\alpha)-\frac12}(\log s)^{h(\alpha)(p(\alpha)-1)}\|f\|_{H^{\alpha}_{\mathsf{mix}}},
\end{align*}
as required.
For $q=\infty$, we apply \cf{l:stech2} with $q=\infty$ and obtain
\begin{align*}
\sigma_s(c)_\infty\lesssim_{\alpha} s^{-h(\alpha)-\frac12}(\log s)^{h(\alpha)(p(\alpha)-1)}\|f\|_{H^{\alpha}_{\mathsf{mix}}},
\end{align*}
which completes the proof.
}

\subsection{Best $s$-term approximation rates in $H^{\beta}$}

We once more introduce some notation that will be used later. Given $\beta=(\beta_1,\beta_2,\dots,\beta_d)\in(0,\infty)^d$ and $r>0$, let
\bes{
C(d, r, \beta) =\left|\left\{n\in\mathbb{Z}^d:\sum_{j=1}^d|n_j|^{\beta_j}<r\right\}\right|.
}

\lem{[]
\label{lem:cdm}
Let $\beta=(\beta_1,\beta_2,\dots,\beta_d)\in(0,\infty)^d$ and $r>0$. Then
\begin{align}\label{beta_res}
C(d,r,\beta)\asymp_{d,\beta}r^\frac{1}{g(\beta)},
\end{align}
where $g(\beta)$ is as in \ef{g-def}. 

}
\prf{
Let \( w = (w_1, \dots, w_d) \in [0,\infty)^d \), $\beta=(\beta_1,\beta_2,\dots,\beta_d)\in(0,\infty)^d$ and $r>0$. Let
\begin{align*}
C_>(d,r,\beta) &= \left| \left\{ n \in \mathbb{Z}_{>}^{d} : \sum_{j=1}^d |n_j|^{\beta_j} < r \right\} \right|,
\\
c(d,r,\beta) &= \mathrm{Vol}_d \left( \left\{ n \in \mathbb{R}^{d} : \sum_{j=1}^d |n_j|^{\beta_j} < r \right\} \right),
\\
c_+(d,r,\beta,w) &= \mathrm{Vol}_d \left( \left\{ n \in \mathbb{R}^{d} : n_j > w_j, \ \forall j \in [d],\   \sum_{j=1}^d |n_j|^{\beta_j} < r \right\} \right),
\end{align*}
where \( \mathrm{Vol}_d(\cdot) \) is the \( d \)-dimensional Lebesgue measure.
Finally, we let $\mathcal{S}(d,s)$ denote the set of all subsets of $[d]$ of size $s$.
Using \cite{wang2005volumes}, we have
\begin{align*}
c(d,1,\beta)= 2^d \frac{\Gamma(1 + 1/\beta_1) \cdots \Gamma(1 + 1/\beta_d)}{\Gamma(1/\beta_1 + 1/\beta_2 + \cdots + 1/\beta_d + 1)},
\end{align*}
and consequently,
\begin{align}\label{wangresult}
c(d,r,\beta)=r^{\frac{1}{g(\beta)}}2^d \frac{\Gamma(1 + 1/\beta_1) \cdots \Gamma(1 + 1/\beta_d)}{\Gamma(1/\beta_1 + 1/\beta_2 + \cdots + 1/\beta_d + 1)}
\end{align}
after a routine calculation.
We now make the following observations:
\begin{enumerate}[label=(\alph*)]
    \item \label{item:pointA}$C(d,r,\beta)=1+\sum_{j=1}^d2^j\left(\sum_{S\in \mathcal{S}(d,j)}C_>(j,r,\beta_S)\right)$, where the 1 accounts for the multi-index $0$.
    \item \label{item:pointB}$c(d,r,\beta)=2^d\times c_+(d,r,\beta,0)$.
    \item \label{item:pointC}$C_>(d,r,\beta)\leq c_+(d,r,\beta,0)$.
    \item \label{item:pointD}$c_+(d,r,\beta,1)\leq C_>(d,r,\beta)$.
    \item \label{item:pointE}$c_+(d,r,\beta,0)\leq 1+\sum_{j=1}^{d-1}\left(\sum_{S\in \mathcal{S}(d,j)}{c_+(j,r,\beta_S,1)}\right)+c_+(d,r,\beta,1)$, where the 1 accounts for the volume of the $d$-dimensional unit cube $[0,1]^d$, and each $c_+(j,r,\beta_S,1)$ (for $S=\left\{s_1,s_2,\dots,s_j\right\}\in\mathcal{S}(d,j)$) accounts for the volume \begin{align*}\mathrm{Vol}_d\left\{n\in\mathbb{R}^d:n_{s_i}>1,\ \forall i \in[j],\ \sum^j_{i=1}|n_{s_i}|^{\beta_{s_{i}}}<r,\ n_s\in[0,1],\  \forall s\notin S\right\}.\end{align*}
\end{enumerate}
We now proceed in two parts, corresponding to the upper bound and lower bound. For the former, we 
combine \ref{item:pointA} and \ref{item:pointC}, followed by \ref{item:pointB}, to get
\begin{align*}
C(d,r,\beta)\leq1+\sum_{j=1}^d2^j\left(\sum_{S\in \mathcal{S}(d,j)}c_+(j,r,\beta_S,0)\right)=1+\sum_{j=1}^d\left(\sum_{S\in \mathcal{S}(d,j)}c(j,r,\beta_S)\right).
\end{align*}
Using \ef{wangresult}, it follows that $C(d,r,\beta)\lesssim_{\beta}r^\frac{1}{g(\beta)}$,
as required.

For the lower bound, we combine \ref{item:pointE} and \ref{item:pointD}, and then apply \ref{item:pointC}, to get
\begin{align*}
c_+(d,r,\beta,0)& \leq 1+\sum_{j=1}^{d-1}\left(\sum_{S\in \mathcal{S}(d,j)}{C_>(j,r,\beta_S)}\right)+C_>(d,r,\beta)
\\
& \leq 1+\sum_{j=1}^{d-1}\left(\sum_{S\in \mathcal{S}(d,j)}{c_+(j,r,\beta_S,0)}\right)+C_>(d,r,\beta)
\end{align*}
Using \ref{item:pointB} we get
\begin{align*}
2^{-d}c(d,r,\beta)\leq 1+\sum_{j=1}^{d-1}\left(\sum_{S\in \mathcal{S}(d,j)}{2^{-j}c(j,r,\beta_S)}\right)+C_>(d,r,\beta)
\end{align*}
With this and \ef{wangresult}, we have $C_>(d,r,\beta)\gtrsim_{\beta}r^\frac{1}{g(\beta)}$.
Finally, we apply \ref{item:pointA} to obtain $C(d,r,\beta)\gtrsim_{\beta}r^\frac{1}{g(\beta)}$,
as required.
}

We may now present a result on best $s$-term approximation in $H^{\beta}$.

\thm{[Best s-term approximation rate in $H^\beta$]\label{t:decay_sum}
Let $\beta=(\beta_1,\beta_2,\dots,\beta_d)\in[0,\infty)^d$ and $g(\beta)$ be as in \ef{g-def}. Let $f\in H^\beta$ and write $c = (\hat{f}_n)_{n \in \bbZ^d}$, where $\hat{f}_n$ is as in \ef{fourier_coefficients_notation}.
Then
\begin{align*}
\sigma_s(c)_q\lesssim_{\beta,q}s^{\frac{1}{q}-g(\beta)-\frac{1}{2}}\nm f_{H^\beta},\quad \forall s \in \bbN,
\end{align*}
for all $q>\frac{2}{2g(\beta)+1}$. Moreover, for $q=\infty$,
\begin{align*}
\sigma_s(c)_\infty\lesssim_{\beta}s^{-g(\beta)-\frac{1}{2}}\nm f_{H^\beta},\quad \forall s \in \bbN.
\end{align*}
}
\prf{
We follow a similar proof to that of \cf{t:decay_prod}. Let $q=2$, $s\in\bbN$, and $\varepsilon>0$ be such that $C(d,1/\varepsilon,\beta)\leq s$, where $C(d,1/\varepsilon,\beta)$ is as in \ef{beta_res}. Then
\be{
\label{eq:tail_bound_sum_sobolev}
\sigma_s(c)^2_2 
\leq \sum_{n : 1 + \sum_{j=1}^d |n_j|^{\beta_j}  \geq \varepsilon^{-1}} | \hat{f}_n |^2  
\leq\varepsilon^2 \sum_{n \in \bbZ^d} \left(1 + \sum_{j=1}^d |n_j|^{\beta_j} \right)^2 | \hat{f}_n |^2 = \varepsilon^2 \|f\|^2_{H^\beta}. 
}
Recall \ef{beta_res}. Lemma \ref{lem:cdm} implies that there are constants $c_1 = c_1(d,\beta),c_2 = c_2(d,\beta) > 0$ such that
\begin{align}
c_1r^\frac{1}{g(\beta)}\leq C(d,r,\beta)\leq c_2r^\frac{1}{g(\beta)}.
\end{align}
Let $\varepsilon$ be such that
\begin{align}
s\leq c_2\left(\frac1\varepsilon-1\right)^\frac{1}{g(\beta)}<s+1\label{eps_btwn_s}
\end{align}
and therefore $C(d,1/\varepsilon,\beta) \leq s$. Notice that $\varepsilon$ satisfies $\varepsilon\leq \left(\frac{s^{g(\beta)}}{c_2}+1\right)^{-1}.
$
Combining this with \ef{eq:tail_bound_sum_sobolev} yields $\sigma_s(c)_2\lesssim_{\beta}s^{-g(\beta)}\|f\|_{H^\beta}$,
which gives the result for $q = 2$.

Next, we apply Stechkin's inequality.
By \cf{l:stech2} with $a=0$, the above bound implies that $c\in w\ell^{\frac{2}{2g(\beta)+1},0}$.
Hence, for $q>\frac{2}{2g(\beta)+1}$, reapplying \cf{l:stech2} with $a=0$ gives
\begin{align*}
\sigma_s(c)_q\lesssim_{\beta,q}s^{\frac{1}{q}-g(\beta)-\frac{1}{2}}\nm f_{H^\beta}.
\end{align*}
For $q=\infty$, we apply \cf{l:stech2} with $a=0$ and $q=\infty$ and obtain
\begin{align*}
\sigma_s(c)_\infty\lesssim_{\beta}s^{-g(\beta)-\frac{1}{2}}\nm f_{H^\beta},
\end{align*}
which completes the proof.
}


\section{Universal algorithms for unknown anisotropy}\label{s:universal-algorithms-thms}

We now construct universal algorithms for recovering functions in the setting of unknown anisotropy, thus establishing contribution (A) of \S \ref{ss:contributions}. 
In this section, $u : \bbN \rightarrow (0,\infty)$ is an arbitrary nondecreasing function with $u(m) \rightarrow \infty$ as $m \rightarrow \infty$. It is informative to think of this function as growing slowly in $m$, e.g., logarithmically ($u(m) = \log(m+1)$) or double-logarithmically ($u(m) = \log(\log(m+1))$).

\subsection{Dominating mixed smoothness spaces $H^{\alpha}_{\mathsf{mix}}$}

Consider the spaces introduced in Definition \ref{def:dom-mixed-space}. Our first two results assert the existence of maps with guaranteed universal approximation properties in probability and in expectation, respectively.

\thm{[Probability bound]
\label{t:main-prob-mixed}
Let $0 < \varepsilon < 1$, $m \geq 2$ and $x_1,\ldots,x_m$ be drawn i.i.d. from the uniform measure on $\bbT^d$. Then there are constants $C(\alpha,d,u)$, $\forall \alpha > 1/2$, and a 
reconstruction map $R : (\bbT^d \times \bbC)^{m} \rightarrow L^2(\bbT^d)$ depending on $u$ and $\varepsilon$ only such that, if
\bes{
\widetilde{m} = \widetilde{m}(\varepsilon,u) = : 
\frac{m}{\log^3(m) u(m) + \log(1/\varepsilon)} \geq 2, 
}
then
\bes{
 \nm{f - R( (x_i,f(x_i))^{m}_{i=1}) }_{L^2}  \leq C(\alpha,d,u) \left ( \frac{\log^{p(\alpha)-1}(\widetilde{m})}{\widetilde{m}} \right )^{h(\alpha)} \nm{f}_{H^{\alpha}_{\mathsf{mix}}}, \quad \forall f \in H^{\alpha}_{\mathsf{mix}}(\bbT^d),  \alpha > 1/2 
}
with probability at least $1-\varepsilon$, where $h$ and $p$ are as in \ef{h-p-def}.
}

\thm{[Expectation bound]\label{t:main-expec-mix}
Let $m \geq 2$ and $x_1,\ldots,x_m$ be drawn i.i.d. from the uniform measure on $\bbT^d$.
Then there are constants $C(\alpha,d,u)$, $\forall \alpha > 1/2$, and a 
reconstruction map 
$
R : (\bbT^d \times \bbC)^{m} \rightarrow L^2(\bbT^d)
$
depending on $u$ only such that, for any bounded set 
$\cA \subset (1/2,\infty)^d$, there exists a constant $C(\cA)$ with the property that, if
\bes{
\widetilde{m} = \widetilde{m}(u) = 
\frac{m}{\log^3(m) u(m) + \log^2(m)} \geq 2, 
}
then
\bes{
\bbE \left [ \sup_{\alpha \in \cA} \sup_{\substack{\substack{f \in H^{\alpha}_{\mathsf{mix}}  \\ f \neq 0} }}  
\frac{\nm{f - R( (x_i,f(x_i))^{m}_{i=1}) }_{L^2} }
{C(\alpha,d,u) \left ( \frac{\log^{p(\alpha)-1}(\widetilde{m})}{\widetilde{m}} \right )^{h(\alpha)} 
\nm{f}_{H^{\alpha}_{\mathsf{mix}}}   }
 \right ] \leq C(\cA),
}
where $h$ and $p$ are as in \ef{h-p-def}.
}

These results establish the existence of a reconstruction map that achieves an approximation rate over $H^{\alpha}_{\mathsf{mix}}$ depending on $h(\alpha)$ and $p(\alpha)$ defined in \ef{h-p-def} and the quantity $\widetilde{m}$. Notice that $\widetilde{m}$ is given by $m$ scaled by a term that can be made arbitrarily close (by choosing $u$) to $\log^3(m)$. Later, in \S \ref{s:lower-bounds-widths}, we will show that this rate is optimal up to this polylogarithmic factor.

The map $R$ asserted by Theorem \ref{t:main-prob-mixed} is universal, since it yields, with high probability, the stipulated rate for \textit{any} value of the anisotropy parameter $\alpha > 1/2$ (see Remark \ref{rem:alpha-half} below for why this condition is needed). The slowly-growing function $u(m)$ is critical in guaranteeing this property.
The map in Theorem \ref{t:main-expec-mix} is universal, but in a slightly weaker sense, since the stipulated rate is only guaranteed for $\alpha$ belonging to some bounded set $\cA$. This arises for technical reasons when the probability bound of Theorem \ref{t:main-prob-mixed} is used to establish the expectation bound in Theorem \ref{t:main-expec-mix}. However, note that $\cA$ can be arbitrary and the map $R$ is independent of $\cA$. 

\rem{
\label{rem:alpha-half}
The condition $\alpha > 1/2$ ensures a continuous embedding $C(\bbT^d) \hookrightarrow H^{\alpha}_{\mathsf{mix}}(\bbT^d)$ and, consequently, the argument of $R$ is well defined for any $f$.
Indeed, for any $f \in H^{\alpha}_{\mathsf{mix}}(\bbT^d)$, we have
\bes{
\|f\|_{L^\infty(\bbT^d)} \le \sum_{n\in\bbZ^d}|\hat f_n|
\le
\Big(\sum_{n\in\bbZ^d}\prod_{j=1}^d(1+|n_j|)^{-2\alpha_j}\Big)^{1/2}
\|f\|_{H^{\alpha}_{\mathsf{mix}}}.
}
The claimed embedding now follows from the fact that
\bes{
\sum_{n\in\bbZ^d}\prod_{j=1}^d(1+|n_j|)^{-2\alpha_j}
=
\prod_{j=1}^d\sum_{k\in\bbZ}(1+|k|)^{-2\alpha_j}<\infty
\quad \Leftrightarrow \quad
\alpha_j>1/2,\  \forall j\in [d] .
}
}

\rem{
\label{rem:log-term}
As we shall see, Theorem \ref{t:main-prob-mixed} (also Theorem \ref{t:main-expec-mix}) yields a rate that is only possibly suboptimal in the factor $\log^3(m) u(m) + \log(1/\varepsilon)$. The term $u(m)$, as discussed above, arises from the requirement that the algorithm be universal and can be chosen to grow arbitrarily slowly. As shown in the proof, the $\log^3(m)$ term arises from requiring the Restricted Isometry Property (RIP) for a certain Fourier matrix of size $m \times N$, where $\log(N) = \ord{\log(m)}$. Reducing the logarithmic term in the RIP for a Fourier matrix is a long-standing open problem \cite{brugiapaglia2021sparse}. Any future reduction would immediately yield a better logarithmic term in our bounds.
}

While our main focus in this paper is on i.i.d.\ samples, in the next result we show the existence of a set of samples that yields truly universal recovery (with no restrictions on $\alpha$) in a non-probabilistic sense. While the rate is unchanged over that of i.i.d.\ sampling, the use of these sample points notably strengthens Theorem \ref{t:main-expec-mix} by removing the need to consider a bounded set $\cA$.

\thm{[Non-i.i.d. sampling]\label{t:non-iid-mix}
Let $m \geq 2$. Then there exists a set of points $(x^{\star}_1,\ldots,x^{\star}_m)$ such that the following holds. There are constants $C(\alpha,d,u)$, $\forall \alpha > 1/2$, and a reconstruction map $R : (\bbT^d \times \bbC)^m \rightarrow L^2(\bbT^d)$ depending on $u$ only such that, if
\bes{
\widetilde{m} = \widetilde{m}(u) = \frac{m}{\log^3(m) u(m)+\log(2)} \geq 2,
}
then
\bes{
 \sup_{\alpha > 1/2} \sup_{\substack{f \in H^{\alpha}_{\mathsf{mix}} \\ f \neq 0}}  \frac{\nm{f - R( (x^{\star}_i,f(x^{\star}_i))^{m}_{i=1}) }_{L^2} }
{C(\alpha,d,u) \left ( \frac{\log^{p(\alpha)-1}(\widetilde{m})}{\widetilde{m}} \right )^{h(\alpha)} \nm{f}_{H^{\alpha}_{\mathsf{mix}}}   }
 \leq 1,
}
where $h$ and $p$ are as in \ef{h-p-def}.
}

\subsection{Anisotropic spaces $H^{\beta}$}

We now present analogous versions of Theorems \ref{t:main-prob-mixed}--\ref{t:non-iid-mix} for the spaces introduced in Definition \ref{def:anisotropic-space}.

\thm{[Probability bound]
\label{t:main-prob-sum}
Let $0<\varepsilon<1$, $m\ge2$ and $x_1,\ldots,x_m$ be drawn i.i.d. from the uniform measure on $\bbT^d$.
Then there exist constants $C(\beta,d,u)$ for all $\beta$ with $g(\beta)>1/2$, where $g$ is as in \ef{g-def}, and a reconstruction map
$R:(\bbT^d\times\bbC)^m\to L^2(\bbT^d)$,
depending only on $u$ and $\varepsilon$, such that if
\[
\widetilde m=\widetilde m(\varepsilon,u)
=\frac{m}{\log^3(m)u(m)+\log(1/\varepsilon)}\ge2,
\]
then
\bes{
\|f-R((x_i,f(x_i))_{i=1}^m)\|_{L^2}
\le
C(\beta,d,u) \widetilde m^{-g(\beta)}\|f\|_{H^\beta},
\forall f\in H^\beta(\bbT^d), \beta : g(\beta)>1/2,
}
with probability at least $1-\varepsilon$.
}

\thm{[Expectation bound]\label{t:main-expec-sum}
Let $m \geq 2$ and $x_1,\ldots,x_m$ be drawn i.i.d. from the uniform measure on $\bbT^d$.
Then there are constants $C(\beta,d,u)$ for all $\beta$ with $g(\beta) > 1/2$, where $g$ is as in \ef{g-def}, and a 
reconstruction map 
$
R : (\bbT^d \times \bbC)^{m} \rightarrow L^2(\bbT^d)
$
depending on $u$ only such that, for any bounded set 
$\cB \subset \left\{\beta\in(1/2,\infty)^d:g(\beta)>1/2\right\}$, there exists a constant $C(\cB)$ with the property that, if
\bes{
\widetilde{m} = \widetilde{m}(u) = 
\frac{m}{\log^3(m) u(m) + \log^2(m) } \geq 2, 
}
then
\bes{
\bbE \left [ \sup_{\beta \in \cB} \sup_{\substack{f \in H^{\beta} \\ f \neq 0}}  
\frac{\nm{f - R( (x_i,f(x_i))^{m}_{i=1}) }_{L^2} }
{C(\beta,d,u) \left ( \frac{1}{\widetilde{m}} \right )^{g(\beta)} 
\nm{f}_{H^{\beta}}   }
 \right ] \leq C(\cB).
}
}
\thm{[Non-i.i.d. sampling]\label{t:non-iid-sum}
Let $m\geq 2$. Then there exists a set of points $(x_1^{\star},\ldots,x_m^{\star})$ such that the following holds. There are constants $C(\beta,d,u)$, for all $\beta$ with $g(\beta)>1/2$, where $g$ is as in \ef{g-def}, and a reconstruction map
$
R:(\bbT^d\times\bbC)^m\to L^2(\bbT^d)
$
depending on $u$ only such that, if
\bes{
\widetilde{m}=\widetilde{m}(u)=\frac{m}{\log^3(m)u(m)+\log(2)}\geq 2,
}
then
\bes{
\sup_{\substack{\beta\in(1/2,\infty)^d\\ g(\beta)>1/2}}
\sup_{\substack{f\in H^{\beta}\\ f\neq 0}}
\frac{\nm{f-R((x_i^{\star},f(x_i^{\star}))_{i=1}^m)}_{L^2}}
{C(\beta,d,u)\Big(\frac{1}{\widetilde{m}}\Big)^{g(\beta)}\nm{f}_{H^{\beta}}}
\leq 1.
}
}

As in the previous case, the rates exhibited in these theorems are optimal up to the distinction between $m$ and $\widetilde{m}$, which is once more polylogarithmic (for a suitable choice of $u$).

\rem{
Just as in the mixed smoothness case, the condition $g(\beta)>1/2$ guarantees the continuous embedding
$H^\beta(\bbT^d)\hookrightarrow C(\bbT^d)$. For succinctness, we omit the details.
}


\subsection{Reduction to a sparse recovery problem}\label{subsec:sparse}

We now turn our attention to the proofs of Theorems \ref{t:main-prob-mixed}--\ref{t:non-iid-mix},\ref{t:main-prob-sum}--\ref{t:non-iid-sum}, commencing with the construction of the algorithms. The main results in the previous section, Theorems \ref{t:decay_prod} and \ref{t:decay_sum}, show that the best $s$-term approximation achieves the desired convergence rates in terms of $s$. Motivated by this, we pursue an approach based on compressed sensing, where we first convert the problem of approximating an unknown $f$ from its samples \ef{f-samples} to a finite-dimensional sparse recovery problem.

Given an unknown function \( f \) and sample points \( x_1, \dots, x_m \) drawn i.i.d. from the uniform distribution on \( \mathbb{T}^d \), we define the normalized sample vector  
\begin{align}\label{eq:noisy_sample_vector_notation}
    b = \frac{1}{\sqrt{m}}\big(f(x_i)\big)_{i=1}^m.
\end{align}
We seek to approximate $f$ by recovering its Fourier coefficients $\hat{f}_\Lambda = (\hat{f}_n)_{n \in \Lambda}$ corresponding to an index set $\Lambda\subseteq\mathbb{Z}^d$ that will be chosen later. Write $N = |\Lambda|$ and let $A\in \mathbb{C}^{m\times N}$ denote the matrix formed by evaluating the truncated Fourier basis $(\phi_n)_{n\in\Lambda}$ at the sample points, divided by $\sqrt m$ to normalize:
\begin{align}\label{eq:Fourier_basis_sampled_randomly}
A = \frac{1}{\sqrt{m}} \left( \phi_n(x_i) \right)_{i \in [m],  n \in \Lambda}.
\end{align}
In addition, we define the truncation error \( v \in \mathbb{C}^m \) as
\begin{align}\label{eq:truncation_error_notation}
v = \frac{1}{\sqrt{m}} \left( f(x_i) - \sum_{n \in \Lambda} \hat{f}_n \phi_n(x_i) \right)_{i \in [m]}.
\end{align}
Substituting $\hat{f}_{\Lambda}$, $A$ and $v$ into \ef{eq:noisy_sample_vector_notation}, we arrive at the following:
\begin{align}\label{eq:pr}
b = A\hat{f}_{\Lambda} + v.
\end{align}
Since $\hat{f}_{\Lambda}$ is approximately sparse, this is a sparse recovery problem. Following \cite{adcock2022sparse,adcock2024efficient,adcock2025optimal,moeller2026best,moeller2025instance,moeller2025sampling}, we solve it using the SR-LASSO decoder \cite{belloni2011square-root}.

\defn{[SR-LASSO]\label{d:srlasso}
Given a regularization parameter $\lambda > 0$, a matrix $A \in \mathbb{C}^{m \times N}$ and vector $b \in \mathbb{C}^m$, the \textit{(unconstrained) Square Root (SR)-LASSO problem} is the optimization problem
\begin{align}\label{eq:srlasso}
\min_{z \in \mathbb{C}^N}  \lambda \|z\|_1 + \|A z - b\|_2.
\end{align}
}

\subsection{Compressed sensing tools}

We now present several key compressed sensing tools. See, e.g., \cite[Chpt.\ 6]{adcock2022sparse}.

\defn{[$\ell^2$-rNSP]
Given $s\in \mathbb{N}$, a matrix $A\in\mathbb{C}^{m\times N}$ satisfies the \textit{$\ell^2$-robust Null Space Property ($\ell^2$-rNSP) of order $s$ with constants $\rho \in (0,1)$ and $\tau > 0$}
if, for all $z\in \mathbb{C}^N$ and all index sets $S\subseteq[N]$ with $|S|=s$,
\begin{align}\label{l2r}
\|z_S\|_2 \leq \frac{\rho}{\sqrt{s}} \|z_{S^c}\|_1 + \tau \|A z\|_2.
\end{align}
}

\lem{[rNSP implies stable and accurate recovery for the SR-LASSO problem]\label{l:srlasso_implies_stable_and_accurate_recovery}
Let $A \in \mathbb{C}^{m \times N}$ satisfy the rNSP of order $s$ with constants $0 < \rho < 1$ and $\tau > 0$. 
Let $z^* \in \mathbb{C}^N$, $h \in \mathbb{C}^m$, and $b = A z^* + h \in \mathbb{C}^m$. 
Then, for any $\lambda$ satisfying $\lambda \in \left(0,\frac{D}{\sqrt s}\right]$ where
$
D = \frac{(1+\rho)}{(3+\rho)\tau}$,
and any $z^\sharp$ such that
\bes{
z^\sharp \in \argmin{z \in \mathbb{C}^N}~\lambda \|z\|_1 + \|A z - b\|_2,
}
we have
\begin{align*}
\|z^* - z^\sharp\|_1 
&\le C_1 \sigma_s(z^*)_1 
+ \frac{1}{2}\left( \frac{C_1}{\lambda} + C_2\sqrt{s} \right)\|h\|_2,
\\
\|z^* - z^\sharp\|_2 
&\le C_3 \frac{\sigma_s(z^*)_1}{\sqrt{s}} 
+ \frac{1}{2}\left( \frac{C_3}{\sqrt{s}\lambda} + C_4 \right)\|h\|_2,
\end{align*}
where the constants $C_1, C_2, C_3, C_4$ are given by
\bes{
C_1 = 2\left( \frac{1+\rho}{1-\rho} \right), \quad
C_2 = \frac{4 \tau}{1-\rho}, \quad
C_3 = \frac{2(1 + \rho)^2}{(1 - \rho)}, \quad
C_4 = 2\tau \frac{(3 + \rho)}{(1 - \rho)}.
}
}

\lem{[The Fourier matrix satisfies the $\ell^2$-rNSP]\label{l:fourier_matrix_satisfies_l2rnsp}
Let $0 < \varepsilon < 1$, $x_1,\ldots,x_m$ be drawn i.i.d. from the uniform distribution on $\bbT^d$ and $A$ be as in \ef{eq:Fourier_basis_sampled_randomly}. Then there is numerical constant $c > 0$ such that, if
\be{
\label{cond-for-rNSP}
m \geq c \cdot s \cdot \left ( \log^2(2s) \cdot \log(2N) + \log(1/\varepsilon) \right )
}
then the matrix $\widetilde{A} = (2 \pi)^{d/2} A$ satisfies the $\ell^2$-rNSP with constants $\rho = \sqrt{2}/3$ and $\tau = 2 \sqrt{5}/3$ with probability at least $1-\varepsilon$.
}
\prf{
The matrix $\widetilde{A}$ is the matrix of a so-called bounded orthonormal system \cite[Def. 6.14]{adcock2022sparse} with constant $K = 1$. Using \cite[Theorem 2.3]{brugiapaglia2021sparse} with $c_{\psi} = C_{\psi} = 1$, we see that the condition
\be{
\label{cond-for-RIP}
m \geq c_1 \cdot s \cdot \log^2(2s) \cdot \log(2N) 
}
implies that $\widetilde{A}$ satisfies the Restricted Isometry Property (RIP) \cite[Definition 6.10]{adcock2022sparse} with constant $\delta = 1/4$, with probability at least $1-2 \exp(-c_2 m / s)$, where $c_1,c_2>0$ are numerical constants. Notice that \ef{cond-for-RIP} holds, due to \ef{cond-for-rNSP}, provided $c$ is sufficiently large. Moreover, \ef{cond-for-rNSP} also implies that $1-2 \exp(-c_2 m / s) \geq 1- \varepsilon$, for sufficiently large $c$. We conclude that $\widetilde{A}$ has the RIP with probability at least $1-\varepsilon$ and constant $\delta = 1/4$. Finally, a standard result \cite[Thm. 6.11]{adcock2022sparse} now implies that $\widetilde{A}$ has the rNSP with constants $\rho = \sqrt{2}/3$ and $\tau = 2 \sqrt{5}/3$, as required.
}

\subsection{Recovery of $H^{\alpha}_{\mathsf{mix}}$ and $H^{\beta}$ functions via the SR-LASSO decoder}

We now establish two theorems that provide error bounds for the recovery of $H^{\alpha}_{\mathsf{mix}}$ and $H^{\beta}$ functions, respectively, from sample values via the SR-LASSO decoder. At this stage, we also specify the truncation set $\Lambda$. For reasons that will become clear later, we now let $\Lambda = \Lambda^{r}_{\mathsf{HC}}$, where
\bes{
\Lambda^{r}_{\mathsf{HC}} = \left\{ n = (n_1,\dots,n_d) \in \mathbb{Z}^d :
\prod_{j=1}^d (1 + |n_j|) \le r \right\}
}
is the hyperbolic cross index set of \textit{order} $r \geq 0$.

\thm{[Error bound for the reconstruction map for anisotropic mixed smoothness Sobolev spaces]\label{t:Error_bound_for_the_reconstruction_map_for_product_type_Sobolev_spaces}
Let $s \in \bbN$ with $s \geq 2$, $r \geq 1$, $0 < \varepsilon < 1$, $\Lambda = \Lambda^{r}_{\mathsf{HC}}$, $N = |\Lambda|$ and consider $x_1,\ldots,x_m$ drawn i.i.d. from the uniform distribution on $\bbT^d$, where $m$ satisfies
\be{
\label{m-cond-rec-map-bd}
m \geq c \cdot s \cdot \left ( \log^2(2s) \cdot \left ( \log(2r) + (d-1) \log(\log(\E r)) \right ) + \log(1/\varepsilon) \right ) 
}
for some universal constant $c > 0$.
Then the following holds with probability at least $1-\varepsilon$. 
Let $\alpha = (\alpha_1, \ldots, \alpha_d) \in (0,\infty)^d$ and $h(\alpha)$ and $p(\alpha)$ be as in \ef{h-p-def}, where $h(\alpha) > 1/2$, $f\in H^{\alpha}_{\mathsf{mix}}$ and suppose that $A$ and $b$ are as in \eqref{eq:noisy_sample_vector_notation} and \eqref{eq:Fourier_basis_sampled_randomly}, respectively. Define $\widetilde{A}=(2\pi)^{d/2}A$ and suppose that
\be{
\label{eq:lambda_dependence_on_s}
 \lambda \in \left(0,\frac{3(\sqrt{2} + 3)}{2\sqrt{5}(\sqrt{2} + 9)\sqrt s}\right] .
}
Then the approximation $f^{\sharp}$ given by
\bes{
f^{\sharp} = \sum_{n\in\Lambda}(2\pi)^{d/2}z^\sharp_n\phi_n, \quad \text{where } z^\sharp\in\argmin{z\in\mathbb{C}^{N}}  ~\lambda\|z\|_1+\|\widetilde Az-b \|_2
}
is well-defined and satisfies
\bes{
\|f - f^\sharp\|_{L^2} \lesssim_{d,\alpha}
        s^{-h(\alpha)} (\log s)^{h(\alpha)(p(\alpha)-1)} 
        \|f\|_{H^{\alpha}_{\mathsf{mix}}}
+ \left(\frac{1}{\sqrt{s}\lambda}+1\right)\left (r^{-h(\alpha)+1/2}\log(r)^{\frac{d-1}{2}}\|f\|_{H^{\alpha}_{\mathsf{mix}}} \right ).
}
}

\prf{
\cf{l:fourier_matrix_satisfies_l2rnsp} states that if
\begin{align}\label{eq:m_grows_as_s_times_polylog_s_times_polylog_N}
m \geq c \cdot s \cdot \left ( \log^2(2s) \cdot \log(2N) + \log(1/\varepsilon) \right ),
\end{align}
then, with probability at least \(1 - \varepsilon\), the matrix \(\widetilde{A}\) satisfies the \(\ell^2\)-rNSP of order \(s\) with constants \(\rho = \tfrac{\sqrt{2}}{3}\) and \(\tau = \tfrac{2\sqrt{5}}{3}\). The hyperbolic cross index set satisfies the bound
\bes{
N = |\Lambda^{\mathsf{HC}}_r| \leq r (1+\log(r))^{d-1} \leq r \log^{d-1}(\E r)
}
for any $r \in \bbN$ and $d \in \bbN$ (see, e.g., \cite[Prop. A.1]{migliorati2013polynomial}). Therefore, the right-hand side of \ef{eq:m_grows_as_s_times_polylog_s_times_polylog_N} satisfies
\bes{
c \cdot s \cdot \left ( \log^2(2s) \cdot \log(2N) + \log(1/\varepsilon) \right ) \leq c \cdot s \cdot \left ( \log^2(2s) \cdot \left ( \log(2r) + (d-1) \log(\log(\E r)) \right ) + \log(1/\varepsilon) \right ) .
}
Hence \ef{eq:m_grows_as_s_times_polylog_s_times_polylog_N} is implied by \ef{m-cond-rec-map-bd}. We deduce that $\widetilde{A}$ has the desired rNSP with probability $1-\varepsilon$.

Now let $\alpha \in (0,\infty)^d$ with $h(\alpha) > 1/2$ and $f \in H^{\alpha}_{\mathsf{mix}}$ be given. Let $v$ be as in \ef{eq:truncation_error_notation}, \( c = (\hat{f}_n)_{n \in \mathbb{Z}^d} \), \( z =(2\pi)^{-d/2} (\hat{f}_n)_{n \in \Lambda} \) and let $c_\Lambda\in\mathbb{C}^{\bbZ^d}$ have the $n$th entry equal to $\hat{f}_n$ if $n \in \Lambda$ and zero otherwise.  The approximation $f^{\sharp}$ is well-defined, since the optimization problem for $z^{\sharp}$ always has a minimizer.
Applying \cf{l:srlasso_implies_stable_and_accurate_recovery} and recalling that $\lambda$ satisfies \ef{eq:lambda_dependence_on_s}, we get that
\be{\label{eq:main_theorem_equation_1}
\begin{split}
\|z - z^\sharp\|_2 
\lesssim  \frac{\sigma_s(z)_1}{\sqrt{s}} + \left(
    \frac{1}{\sqrt{s}\lambda}
    + 1
  \right)\|v\|_2 .
  \end{split}
}
Additionally, note that
\bes{
\|f - f^\sharp\|_{L^2} \le (2\pi)^{d/2}\|z - z^\sharp\|_2 + \|f - \sum_{n\in\Lambda}\hat f_n\phi_n\|_{L^2} \leq (2\pi)^{d/2}\|z - z^\sharp\|_2 + \nm{c-c_{\Lambda}}_1
}
and
\begin{align*}
\nm{v}_2 \leq \| f - \sum_{j \in \Lambda} \hat{f}_j \phi_j \|_{L^{\infty}} \leq \nm{c-c_{\Lambda}}_1 .
\end{align*}
Combining this with the previous expression we deduce that
\be{
\label{f-fsharp-s-Lambda}
\|f - f^\sharp\|_{L^2} \lesssim_{d}  \frac{\sigma_s(z)_1}{\sqrt{s}} + \left(
    \frac{1}{\sqrt{s}\lambda}
    + 1
  \right) \nm{c-c_{\Lambda}}_1 .
}
To conclude, we bound the two terms on the right-hand side. For the first term, we invoke \cf{t:decay_prod} with $q=1$. Since $h(\alpha)>1/2$, we have
\be{
\label{sigma-s-z-bd}
\frac{\sigma_s(z)_1}{\sqrt{s}} = (2\pi)^{-d/2}\frac{\sigma_s(c)_1}{\sqrt{s}}\lesssim_{d,\alpha} s^{-h(\alpha)}(\log s)^{h(\alpha)(p(\alpha)-1)}\|f\|_{H^{\alpha}_{\mathsf{mix}}} .
}
For the second term of \ef{f-fsharp-s-Lambda}, we proceed as follows. Writing
\[
\|c-c_{\Lambda}\|_{1}=\sum_{n\notin\Lambda}|c_{n}| = \sum_{n\notin\Lambda}
|c_{n}| 
\frac{\prod_{j=1}^{d}\max\{1,|n_{j}|\}^{h(\alpha)}}{\prod_{j=1}^{d}\max\{1,|n_{j}|\}^{h(\alpha)}}.
\]
and applying the Cauchy--Schwarz inequality gives
\begin{align}\label{eq:an_inequality_in_product_type_proof}
\|c-c_{\Lambda}\|_{1}
\le 
\Biggl(\sum_{n\notin\Lambda} |c_{n}|^{2}\prod_{j=1}^{d}\max\{1,|n_{j}|\}^{2h(\alpha)}\Biggr)^{1/2}
\Biggl(\sum_{n\notin\Lambda}\prod_{j=1}^{d}\max\{1,|n_{j}|\}^{-2h(\alpha)}\Biggr)^{1/2}.
\end{align}
For the first factor, note that $h(\alpha)\le \alpha_j$ for each $j\in[d]$ and that 
$\max\{1,|n_j|\}\le 1+|n_j|$. 
Hence
\[
\sum_{n\notin\Lambda} |c_{n}|^{2}\prod_{j=1}^{d}\max\{1,|n_{j}|\}^{2h(\alpha)}
\le \sum_{n\in\bbZ^d} |c_{n}|^{2}\prod_{j=1}^{d}(1+|n_{j}|)^{2\alpha_j}
= \|f\|_{H^{\alpha}_{\mathsf{mix}}}^{2}.
\]
For the second factor, we have $\left\{n\in\bbZ^{d}: \prod_{j=1}^{d}\max\{1,|n_{j}|\}\le \frac{r}{2^{d}}\right\}
\subseteq \Lambda$.
Thus,
\[
\sum_{n\notin\Lambda}\prod_{j=1}^{d}\max\{1,|n_{j}|\}^{-2h(\alpha)}
\le 
\sum_{\substack{n\in\bbZ^{d}: \prod_{j=1}^{d}\max\{1,|n_{j}|\}> \frac{r}{2^{d}}}}
\prod_{j=1}^{d}\max\{1,|n_{j}|\}^{-2h(\alpha)}.
\]
By \cite[Theorem~2.30]{adcock2010modified}, for $2h(\alpha)>1$,
\[
\sum_{\substack{n\in\bbZ^{d}: \prod_{j=1}^{d}\max\{1,|n_{j}|\}> \frac{r}{2^{d}}}}
\prod_{j=1}^{d}\max\{1,|n_{j}|\}^{-2h(\alpha)}
\lesssim_{d,\alpha} r^{-2h(\alpha)+1}\log(r)^{d-1} .
\]
Combining the above estimates gives $
\|c-c_{\Lambda}\|_{1}
\lesssim_{d,\alpha} r^{-h(\alpha)+\frac12}\log(r)^{\frac{d-1}{2}}
\|f\|_{H^{\alpha}_{\mathsf{mix}}} .
$
Substituting this and \ef{sigma-s-z-bd} into \ef{f-fsharp-s-Lambda} now completes the proof.
}

\thm{[Error bound for the reconstruction map for anisotropic Sobolev spaces]
\label{t:Error_bound_for_the_reconstruction_map_for_sum_type_Sobolev_spaces}
Let $s \in \bbN$ with $s \geq 2$, $r \geq 1$, $0 < \varepsilon < 1$, $\Lambda = \Lambda^{r}_{\mathsf{HC}}$, $N = |\Lambda|$ and consider $x_1,\ldots,x_m$ drawn i.i.d. from the uniform distribution on $\bbT^d$, where $m$ satisfies
\begin{align*}
m \geq
 c s \left ( 
\log^{2}(2s) \log r
+ (d - 1) \log^{2}(2s) \log\log r 
+ \log^{2}(2s) \log\frac{2}{(d - 1)!}
+ \log^{2}(2s)
+ \log\frac{1}{\varepsilon}
\right )
\end{align*}
for some absolute constant $c > 0$.
Then the following holds with probability at least $1-\varepsilon$.
Let $\beta = (\beta_1, \ldots, \beta_d) \in (0,\infty)^d$ and
$g(\beta) > 1/2$ be as in \ef{g-def}, $f\in H^\beta$ and suppose that $A$ and $b$ are as in
\eqref{eq:noisy_sample_vector_notation} and \eqref{eq:Fourier_basis_sampled_randomly}, respectively.
Define $\widetilde{A}=(2\pi)^{d/2}A$ and suppose that
\be{
\label{eq:lambda_dependence_on_s_sum_type}
 \lambda \in \left(0,\frac{3(\sqrt{2} + 3)}{2\sqrt{5}(\sqrt{2} + 9)\sqrt s}\right] .
}
Then the approximation $f^{\sharp}$ given by
\bes{
f^\sharp=\sum_{n\in\Lambda}(2\pi)^{d/2}z^\sharp_n\phi_n,
\quad
z^\sharp\in\argmin{z\in\mathbb{C}^{N}}~  \lambda\|z\|_1+\|\widetilde Az-b\|_2
}
is well-defined and satisfies
\bes{
\|f - f^\sharp\|_{L^2} \lesssim_{d,\beta}
        s^{-g(\beta)} 
        \|f\|_{H^\beta}
+ \left(\frac{1}{\sqrt{s}\lambda}+1\right)
        \left(
        r^{-g(\beta)+1/2}\log(r)^{\frac{d-1}{2}}\|f\|_{H^\beta}
        \right).
}
}

\prf{
The proof follows the same structure as that of \cf{t:Error_bound_for_the_reconstruction_map_for_product_type_Sobolev_spaces}. Recall from \ef{f-fsharp-s-Lambda} that
\be{
\label{f-fsharp-s-Lambda-FOR-BETA-SPACES}
\|f - f^\sharp\|_{L^2} \lesssim_{d}  \frac{\sigma_s(z)_1}{\sqrt{s}} + \left(
    \frac{1}{\sqrt{s}\lambda}
    + 1
  \right) \nm{c-c_{\Lambda}}_1 .
}
with probability at least $1-\varepsilon$ for $\lambda$ satisfying \ef{eq:lambda_dependence_on_s} whenever $m$ satisfies \ef{m-cond-rec-map-bd}.
We first estimate the term $\sigma_s(z)_1/\sqrt{s}$ using \cf{t:decay_sum} with $q=1$.
Since $g(\beta)>\tfrac12$, this gives
\be{
\label{sigma-s-z-bd-sum}
\frac{\sigma_s(z)_1}{\sqrt{s}}
= (2\pi)^{-d/2}\frac{\sigma_s(c)_1}{\sqrt{s}}
\lesssim_{d,\beta} s^{-g(\beta)}\|f\|_{H^\beta}.
}
We now estimate the term $\|c-c_{\Lambda}\|_1$.
As in the proof of
\cf{t:Error_bound_for_the_reconstruction_map_for_product_type_Sobolev_spaces},
we write
\[
\|c-c_{\Lambda}\|_{1}
=\sum_{n\notin\Lambda}|c_{n}|
=\sum_{n\notin\Lambda}
|c_{n}|
\frac{\prod_{j=1}^{d}\max\{1,|n_{j}|\}^{g(\beta)}}{\prod_{j=1}^{d}\max\{1,|n_{j}|\}^{g(\beta)}},
\]
and apply the Cauchy--Schwarz inequality to obtain
\begin{align}\label{eq:an_inequality_in_sum_type_proof}
\|c-c_{\Lambda}\|_{1}
\le 
\Biggl(\sum_{n\notin\Lambda} |c_{n}|^{2}\prod_{j=1}^{d}\max\{1,|n_{j}|\}^{2g(\beta)}\Biggr)^{1/2}
\Biggl(\sum_{n\notin\Lambda}\prod_{j=1}^{d}\max\{1,|n_{j}|\}^{-2g(\beta)}\Biggr)^{1/2}.
\end{align}
For the first factor, we use the bound $\max\{1,|n_j|\}\le 1+|n_j|$ to get
\begin{align*}
\sum_{n\notin\Lambda} |c_{n}|^{2}\prod_{j=1}^{d}\max\{1,|n_{j}|\}^{2g(\beta)}
&\le \sum_{n\in\bbZ^d} |c_{n}|^{2}\prod_{j=1}^{d}(1+|n_{j}|)^{2g(\beta)} 
\\
& \lesssim_{d,\beta} \sum_{n\in\bbZ^d} \left( 1 + \sum_{j=1}^{d}|n_j|^{\beta_j} \right)^{2} |c_n|^2 =  \|f\|_{H^\beta}^{2}.
\end{align*} 
Here we used that $\prod_{j=1}^{d}(1+|n_j|)^{g(\beta)} \lesssim_{d,\beta} 1 + \sum_{j=1}^{d}|n_j|^{\beta_j}$ for all $n\in\bbZ^d$, which follows from the weighted AM-GM inequality. Indeed, let $w_j=g(\beta)/\beta_j$ for $j\in[d]$ so that $w_j>0$ and $\sum_{j=1}^{d}w_j=1$. For $n\in\bbZ^d$, we apply the weighted AM-GM inequality with $a_j=(1+|n_j|)^{\beta_j}$ to get
\begin{align*}
\prod_{j=1}^{d}(1+|n_j|)^{g(\beta)} 
=\prod_{j=1}^{d}a_j^{w_j} 
\le \sum_{j=1}^{d} w_j a_j 
=\sum_{j=1}^{d} w_j (1+|n_j|)^{\beta_j}.
\end{align*}
Since $1+|t|\le 2\max\{1,|t|\}$ for $t\ge 0$, we have $(1+|t|)^{\beta}\le 2^{\beta}(1+|t|^{\beta})$ for all $\beta>0$. Therefore,
\begin{align*}
\sum_{j=1}^{d} w_j (1+|n_j|)^{\beta_j} 
\le 2^{\max_{k\in[d]}\beta_k}\sum_{j=1}^{d} w_j (1+|n_j|^{\beta_j})
\le 2^{\max_{k\in[d]}\beta_k}\left(1+\sum_{j=1}^{d}|n_j|^{\beta_j}\right).
\end{align*}
This gives $\prod_{j=1}^{d}(1+|n_j|)^{g(\beta)} \lesssim_{d,\beta} 1+\sum_{j=1}^{d}|n_j|^{\beta_j}$, as claimed.

For the second factor in \ef{eq:an_inequality_in_sum_type_proof}, note that
$
\left\{n\in\bbZ^{d}: \prod_{j=1}^{d}\max\{1,|n_{j}|\}\le \frac{r}{2^{d}}\right\}
\subseteq \Lambda,
$
which implies that
\[
\sum_{n\notin\Lambda}\prod_{j=1}^{d}\max\{1,|n_{j}|\}^{-2g(\beta)}
\le 
\sum_{\substack{n\in\bbZ^{d}: \prod_{j=1}^{d}\max\{1,|n_{j}|\}> \frac{r}{2^{d}}}}
\prod_{j=1}^{d}\max\{1,|n_{j}|\}^{-2g(\beta)}.
\]
By \cite[Theorem~2.30]{adcock2010modified}, for $2g(\beta)>1$,
\[
\sum_{\substack{n\in\bbZ^{d}: \prod_{j=1}^{d}\max\{1,|n_{j}|\}> \frac{r}{2^{d}}}}
\prod_{j=1}^{d}\max\{1,|n_{j}|\}^{-2g(\beta)}
\lesssim_{d,\beta} r^{-2g(\beta)+1}\log(r)^{d-1}.
\]
Combining the above estimates yields
$
\|c-c_{\Lambda}\|_{1}
\lesssim_{d,\beta} r^{-g(\beta)+\frac12}\log(r)^{\frac{d-1}{2}}
\|f\|_{H^{\beta}} .
$
Substituting this and \ef{sigma-s-z-bd-sum} into \ef{f-fsharp-s-Lambda-FOR-BETA-SPACES} now gives the result.
}

\subsection{Proofs of Theorems \ref{t:main-prob-mixed}--\ref{t:non-iid-mix},\ref{t:main-prob-sum}--\ref{t:non-iid-sum}}

We are now ready to prove the main results of this section. In all cases, we build a reconstruction map based on the SR-LASSO decoder, while making judicious choices for the various parameters ($r$, $\lambda$ and so forth) so as to guarantee the desired approximation rates.

\prf{[Proof of Theorem~\ref{t:main-prob-mixed}]

We begin by invoking Theorem~\ref{t:Error_bound_for_the_reconstruction_map_for_product_type_Sobolev_spaces}. 
We introduce positive integers $r = r(m,\varepsilon,u)$ and $s = s(m,\varepsilon,u)$, whose dependence on $m$, $\varepsilon$ and $u$ will be specified later.  Let $\Lambda = \Lambda^r_{\mathsf{HC}}$ and $x_1, \ldots, x_m$ be drawn i.i.d. from the uniform measure on $\mathbb{T}^d$. 
Let $b$ and $A$ be as in \ef{eq:noisy_sample_vector_notation} and \ef{eq:Fourier_basis_sampled_randomly}, respectively, and set $\widetilde{A} = (2\pi)^{d/2} A$ and $\lambda = \frac{3(\sqrt{2} + 3)}{2\sqrt{5}(\sqrt{2} + 9)\sqrt{s}}$.
We then define the map
\[
R\big((x_i, f(x_i))_{i=1}^m\big) = f^\sharp = \sum_{n \in \Lambda} (2\pi)^{d/2} z_n^\sharp \phi_n,
\]
where $z^\sharp$ is the minimal $2$-norm solution of the SR-LASSO problem based on $\widetilde{A}$, $b$ and $\lambda$: namely,
\[
z^\sharp = \argmin{} \left \{ \nm{z^*}_2 : z^* \in \argmin{z \in \mathbb{C}^N}~\lambda \|z\|_1 + \|\widetilde{A}z - b\|_2 \right \}.
\]
Note that $z^{\sharp}$ exists (since the SR-LASSO problem has a minimizer) and is unique, since the set of minimizers is a convex set and $\nm{\cdot}_2$ is strictly convex. Hence $R$ is well defined.

Theorem \ref{t:Error_bound_for_the_reconstruction_map_for_product_type_Sobolev_spaces} and the above choice of $\lambda$ implies that if
\be{
\label{m-cond-main-prob-bd}
m \geq c \cdot s \cdot \left ( \log^2(2s) \cdot \left ( \log(2r) + (d-1) \log(\log(\E r)) \right ) + \log(1/\varepsilon) \right ) 
}
for some universal constant $c > 0$, then 
\be{
\label{eq:error_bound_for_theorem_1point1}
\nm{f - R( (x_i,f(x_i))^{m}_{i=1}) }_{L^2} \lesssim_{d,\alpha} \left ( s^{-h(\alpha)} (\log s)^{h(\alpha)(p(\alpha)-1)} 
      + r^{-h(\alpha)+1/2}\log(r)^{\frac{d-1}{2}} \right ) \|f\|_{H^{\alpha}_{\mathsf{mix}}} .
}
Let $r = \lceil s^{u(m)} \rceil$ and
\eas{
s & = \max \left \{ \left \lfloor \frac{m}{c' \left ( \log^2(2m) \cdot \left (\log(2) +  u(m) \log(m) + (d-1) \log (
1 + u(m) \log(m) ) \right )  + \log(1/\varepsilon) \right ) } \right \rfloor , 1 \right \},
}
where $c' > 0$ is a universal constant which is chosen sufficiently large so that $c' \geq c$ (the constant in \ef{m-cond-main-prob-bd}) and $s \leq m$.

We now consider two cases: $s \geq 2$ and $s = 1$. Suppose first that $s \geq 2$. Then \ef{m-cond-main-prob-bd} holds for this choice of $s$ and $r$. 
Since $s \geq 2$ by assumption, we have
\eas{
r^{-h(\alpha) + 1/2} \log(r)^{\frac{d-1}{2}}  \leq s^{-u(m)(h(\alpha)-1/2)}  \log(2 s^{u(m)})^{\frac{d-1}{2}}
& \leq s^{-u(m)(h(\alpha)-1/2)} (u(m)+1)^{\frac{d-1}{2}} \log(s)^{\frac{d-1}{2}}
\\
&  \leq s^{\frac{d-1}{2} - u(m)(h(\alpha)-1/2) + \frac{d-1}{2} \log (u(m)+1) / \log(2) }
}
and since $u(m) \rightarrow \infty$, we see that there exists an $m = m_1(d,\alpha,u)$ such that
\bes{
r^{-h(\alpha) + 1/2} \log(r)^{\frac{d-1}{2}} \leq s^{-h(\alpha)} (\log(s))^{h(\alpha)(p(\alpha)-1)}
}
for all $m \geq m_1(d,\alpha,u)$. Hence, this and \ef{eq:error_bound_for_theorem_1point1} give that
\bes{
\nm{f - R( (x_i,f(x_i))^{m}_{i=1}) }_{L^2} \lesssim_{d,\alpha} s^{-h(\alpha)} (\log(s))^{h(\alpha)(p(\alpha)-1)}  \|f\|_{H^{\alpha}_{\mathsf{mix}}} 
}
if $m \geq m_1(d,\alpha,u)$. Moreover, when $s \geq 2$, we also have that $s \asymp_d \widetilde{m}$. Hence 
\bes{
\nm{f - R( (x_i,f(x_i))^{m}_{i=1}) }_{L^2} \lesssim_{d,\alpha} \left ( \frac{\log^{p(\alpha)-1}(\widetilde{m})}{\widetilde{m}} \right )^{h(\alpha)}  \|f\|_{H^{\alpha}_{\mathsf{mix}}}.
}
Conversely, if $s \geq 2$ but $m < m_1(d,\alpha,u)$ then, since $r \geq 1$, we have
\bes{
\nm{f - R( (x_i,f(x_i))^{m}_{i=1}) }_{L^2} \lesssim_{d,\alpha} \nm{f}_{H^{\alpha}_{\mathsf{mix}}} 
}
We deduce that, for $s \geq 2$,
\bes{
 \nm{f - R( (x_i,f(x_i))^{m}_{i=1}) }_{L^2} \lesssim_{d,\alpha} \nm{f}_{H^{\alpha}_{\mathsf{mix}}} \lesssim_{d,\alpha} \begin{cases} \nm{f}_{H^{\alpha}_{\mathsf{mix}}} & m < m_1(d,\alpha,u)  \\ \left ( \frac{\log^{p(\alpha)-1}(\widetilde{m})}{\widetilde{m}} \right )^{h(\alpha)} \|f\|_{H^{\alpha}_{\mathsf{mix}}}  & m \geq m_1(d,\alpha,u).  \end{cases} 
}
However, since $\widetilde{m} \geq 2$ by assumption, this immediately implies that 
\bes{
 \nm{f - R( (x_i,f(x_i))^{m}_{i=1}) }_{L^2} \lesssim_{d,\alpha,u} \left ( \frac{\log^{p(\alpha)-1}(\widetilde{m})}{\widetilde{m}} \right )^{h(\alpha)} \|f\|_{H^{\alpha}_{\mathsf{mix}}}.
}
Now suppose that $s = 1$. 
By virtue of the fact that $z^{\sharp}$ is a minimizer of the SR-LASSO problem, we have
\bes{
\nm{f^{\sharp}}_{L^2} \lesssim_{d} \nm{z^{\sharp}}_1 \leq \frac{1}{\lambda} \left ( \lambda \nm{0}_1 + \nm{\widetilde{A} 0 - b }_2 \right ) \lesssim \sqrt{s} \nm{b}_2 \leq \nm{f}_{L^{\infty}} \lesssim_{d,\alpha}  \nm{f}_{H^{\alpha}_{\mathsf{mix}}}.
}
Here, in the penultimate step, we used the continuous embedding $ H^{\alpha}_{\mathsf{mix}}(\bbT^d)\hookrightarrow C(\bbT^d)$. Therefore
\bes{
\nm{f - R( (x_i,f(x_i))^{m}_{i=1}) }_{L^2} \leq \nm{f}_{L^2} + \nm{f^{\sharp}}_{L^2} \lesssim_{d,\alpha} \nm{f}_{H^{\alpha}_{\mathsf{mix}}}.
}
But, when $s = 1$ we have
\bes{
\frac{m}{c' \cdot \left ( \log^2(2m) \cdot \left (\log(2) +  u(m) \log(m) + (d-1) \log (
1 + u(m) \log(m) ) \right )  + \log(1/\varepsilon) \right ) } < 1
}
and therefore $\widetilde{m} \lesssim_{d} 1$.
Since $\widetilde{m} \geq 2$ by assumption, we deduce that
\bes{
\nm{f - R( (x_i,f(x_i))^{m}_{i=1}) }_{L^2} \lesssim_{d,\alpha} \left( \frac{\log^{p(\alpha)-1}(\widetilde{m})}{\widetilde{m}} \right)^{h(\alpha)} \|f\|_{H^{\alpha}_{\mathsf{mix}}},
}
in this case as well. This completes the proof.
}

\prf{[Proof of Theorem~\ref{t:main-expec-mix}]
Let $\varepsilon = m^{-\log(m)} < 1$ and consider the reconstruction map $R$ from Theorem \ref{t:main-prob-mixed} with associated constants $C(\alpha,d,u)$. We may assume that $C(\alpha,d,u) \geq 1$. Let $X$ denote the random variable
\bes{
X = \sup_{\alpha \in \cA} \sup_{\substack{f \in H^{\alpha}_{\mathsf{mix}} \\ f \neq 0}}  \frac{\nm{f - R( (x_i,f(x_i))^{m}_{i=1}) }_{L^2} }
{C(\alpha,d,u) \left ( \frac{\log^{p(\alpha)-1}(\widetilde{m})}{\widetilde{m}} \right )^{h(\alpha)} \nm{f}_{H^{\alpha}_{\mathsf{mix}}}   }.
}
We want to show that $\bbE[X] \leq 1$.
By Theorem \ref{t:main-prob-mixed}, we have $\bbP(X > 1 ) \leq \varepsilon$.
Hence, by the law of total expectation, we have
\eas{
\bbE \left [ X \right ] &= \bbE \left [ X | X>1 \right ] \bbP(X > 1) + \bbE \left [ X | X \leq 1 \right ] \bbP(X \leq 1) \leq \bbE [X | X > 1 ] \varepsilon + 1.
}
Now, from arguments given in the previous proof, we know that
\bes{
\nm{f - R( (x_i,f(x_i))^{m}_{i=1}) }_{L^2} \lesssim_{\alpha,d} \nm{f}_{H^{\alpha}_{\mathsf{mix}} }.
}
Therefore
\bes{
X \leq  \sup_{\alpha \in \cA} \frac{1}{ C(\alpha,d,u) \left ( \frac{\log^{p(\alpha)-1}(\widetilde{m})}{\widetilde{m}} \right )^{h(\alpha)}} \lesssim  \sup_{\alpha \in \cA} \frac{\widetilde{m}^{h(\alpha)} }{ 2^{h(\alpha)(p(\alpha)-1)} } \leq c_1(\cA)  \cdot \widetilde{m}^{c_2(\cA)} 
}
Using the definition of $\varepsilon$, we deduce that
\bes{
\bbE[X] \leq c_1(\cA) \widetilde{m}^{c_2(\cA)} m^{-\log(m)} + 1 \lesssim c_1(\cA) m^{c_2(\cA) - \log(m)} + 1 \leq C(\cA).
}
This gives the result.
}
\prf{[Proof of Theorem~\ref{t:non-iid-mix}]
Consider the reconstruction map from Theorem \ref{t:main-prob-mixed} with $\varepsilon = 1/2$. We now proceed similarly, working through the two cases $s \geq 2$ and $s =1$ separately.
Suppose first that $s \geq 2$. Then 
\be{
\label{eq:error_bound_for_det_thm}
\nm{f - R( (x_i,f(x_i))^{m}_{i=1}) }_{L^2} \lesssim_{d,\alpha,u} 
\left ( \frac{\log^{p(\alpha)-1}(\widetilde{m})}{\widetilde{m}} \right )^{h(\alpha)} \nm{f}_{H^{\alpha}_{\mathsf{mix}}}
 \quad \forall f \in H^{\alpha}_{\mathsf{mix}}(\bbT^d), \alpha > 1/2
}
with probability at least $1/2$. This implies the existence of points $(x^{\star}_1,\ldots,x^{\star}_m)$ such that \ef{eq:error_bound_for_det_thm} holds for $x_i = x^{\star}_i$.

Now suppose that $s = 1$. In this case, we pick any points $(x^{\star}_1,\ldots,x^{\star}_m) \subset \bbT^d$. The argument given in the proof of Theorem \ref{t:main-prob-mixed} for this case makes no assumption on the sample points and leads to the same bound \ef{eq:error_bound_for_det_thm}. This completes the proof.}

\prf{[Proof of Theorem~\ref{t:main-prob-sum}]
We argue along the same lines as in the proof of
Theorem~\ref{t:main-prob-mixed} and invoke Theorem~\ref{t:Error_bound_for_the_reconstruction_map_for_sum_type_Sobolev_spaces}.
Let $r=r(m,\varepsilon,u)$ and $s=s(m,\varepsilon,u)$ be positive integers to be specified below, and set
$\Lambda=\Lambda_{\mathsf{HC}}^r$.
Let $x_1,\ldots,x_m$ be drawn i.i.d. from the uniform measure on $\bbT^d$.
Let $A$ and $b$ be as in
\eqref{eq:noisy_sample_vector_notation} and \eqref{eq:Fourier_basis_sampled_randomly}, respectively,
and set $\widetilde A=(2\pi)^{d/2}A$ and $
\lambda=\frac{3(\sqrt2+3)}{2\sqrt5(\sqrt2+9)\sqrt s}$. We then define the map
\[
R((x_i,f(x_i))_{i=1}^m)
=f^\sharp=\sum_{n\in\Lambda}(2\pi)^{d/2}z_n^\sharp\phi_n,
\]
where $z^\sharp$ is chosen to be the minimal $2$-norm minimizer of the SR-LASSO problem based on $\widetilde A$, $b$ and $\lambda$, i.e.,
\[
z^\sharp=\arg\min\Bigl\{\|z^*\|_2:z^*\in\argmin{z\in\bbC^N}~\lambda\|z\|_1+\|\widetilde Az-b\|_2\Bigr\}.
\]
As in the proof of Theorem~\ref{t:main-prob-mixed}, this choice is well defined and unique, and hence $R$ is well defined.

By Theorem~\ref{t:Error_bound_for_the_reconstruction_map_for_sum_type_Sobolev_spaces} and the above choice of $\lambda$,
there exists a constant $c>0$ such that if
\be{\label{label-for-constant-from-rNSP-lower-bound-for-m}
m\ge
c s\Bigl(\log^2(2s)(\log(2r)+(d-1)\log(\log(e r)))+\log(1/\varepsilon)\Bigr),
}
then, with probability at least $1-\varepsilon$,
\be{\label{f-fsharp-error-rate-s-and-r-Hbeta-norm}
\|f-f^\sharp\|_{L^2}
\lesssim_{d,\beta}
\Bigl(
s^{-g(\beta)}
+r^{-g(\beta)+1/2}\log(r)^{(d-1)/2}
\Bigr)\|f\|_{H^\beta}.
}
Let $r=\lceil s^{u(m)}\rceil$ and  
\eas{
s=\max\Biggl\{
\Biggl\lfloor
\frac{m}{
c'\Bigl(
\log^2(2m)\bigl(\log2+u(m)\log m+(d-1)\log(1+u(m)\log m)\bigr)
+\log(1/\varepsilon)
\Bigr)}
\Biggr\rfloor,
1
\Biggr\},
}
where $c'>0$ is chosen so that $c'\ge c$, where $c$ is as in
\eqref{label-for-constant-from-rNSP-lower-bound-for-m},
and $s\le m$.
Now, following along the proof of Theorem~\ref{t:main-prob-mixed}, we consider the cases
$s\ge2$ and $s=1$. Suppose first that $s\ge2$. Then \ef{label-for-constant-from-rNSP-lower-bound-for-m}
holds for this choice of $s$ and $r$. Using that $u(m)\to\infty$ and that $s\ge2$, we argue as in the
proof of Theorem~\ref{t:main-prob-mixed} to conclude that there exists $m_1=m_1(d,\beta,u)$ such that $r^{-g(\beta)+1/2}\log(r)^{(d-1)/2}\le s^{-g(\beta)}$
for all $m\ge m_1$. Combining this with \ef{f-fsharp-error-rate-s-and-r-Hbeta-norm}, we obtain
$
\|f-f^\sharp\|_{L^2}\lesssim_{d,\beta}s^{-g(\beta)}\|f\|_{H^\beta}
$
for all $m\ge m_1$.
The remainder of the argument for $s\ge2$ follows along the lines
of the proof of Theorem~\ref{t:main-prob-mixed}, and we ultimately conclude that
\[
\|f-R((x_i,f(x_i))_{i=1}^m)\|_{L^2}
\lesssim_{d,\beta,u}\widetilde m^{-g(\beta)}\|f\|_{H^\beta}
\]
whenever $s\ge 2$.
Now suppose that $s=1$. The same calculation as in the proof of
Theorem~\ref{t:main-prob-mixed} shows that
$
\|f^\sharp\|_{L^2}\lesssim_{d,\beta}\|f\|_{H^\beta}.
$
Moreover, when $s=1$ we have $\widetilde m\lesssim_d 1$.
Since $\widetilde m\ge2$ by assumption, we deduce that
\[
\|f-R((x_i,f(x_i))_{i=1}^m)\|_{L^2}
\lesssim_{d,\beta,u}\widetilde m^{-g(\beta)}\|f\|_{H^\beta}.
\]
This completes the proof.
}

\prf{[Proof of Theorem~\ref{t:main-expec-sum}]
The proof follows along the lines of the proof of Theorem \ref{t:main-expec-mix}, with Theorem \ref{t:main-prob-sum} in place of Theorem \ref{t:main-prob-mixed}.
Let $\varepsilon=m^{-\log(m)}<1$ and take the reconstruction map $R$ from Theorem \ref{t:main-prob-sum} with constants $C(\beta,d,u)$. We may assume $C(\beta,d,u)\geq 1$. Define
\bes{
X=\sup_{\beta\in\cB}\sup_{\substack{f\in H^\beta \\ f\neq 0}}
\frac{\|f-R((x_i,f(x_i))_{i=1}^m)\|_{L^2}}
{C(\beta,d,u)\Big(\frac{1}{\widetilde{m}}\Big)^{g(\beta)}\|f\|_{H^\beta}}.
}
By Theorem \ref{t:main-prob-sum}, $\bbP(X>1)\leq \varepsilon$.
Applying the law of total expectation yields
\bes{
\bbE[X]\leq \bbE[X\mid X>1]\varepsilon+1.
}
Next, as in the proof of Theorem \ref{t:main-expec-mix}, we rely on a crude upper bound. Since $z^{\sharp}$ is a minimizer of the SR-LASSO problem, we can compare it with the zero vector to get
\bes{
\|f^{\sharp}\|_{L^2}\lesssim_{d}\|z^{\sharp}\|_{1}
\leq \frac{1}{\lambda}\Big(\lambda\|0\|_{1}+\|\widetilde{A}0-b\|_{2}\Big)
\lesssim \sqrt{s}\,\|b\|_{2}
\leq \|f\|_{L^{\infty}}
\lesssim_{d,\beta}\|f\|_{H^{\beta}}.
}
Here, in the penultimate step, we used the continuous embedding $H^{\beta}(\bbT^d)\hookrightarrow C(\bbT^d)$. Consequently,
\bes{
\|f-R((x_i,f(x_i))_{i=1}^{m})\|_{L^2}
\leq \|f\|_{L^2}+\|f^{\sharp}\|_{L^2}
\lesssim_{d,\beta}\|f\|_{H^{\beta}}.
}
Therefore,
\bes{
X\le \sup_{\beta\in\cB}\frac{1}{C(\beta,d,u)}\widetilde{m}^{g(\beta)}
\lesssim c_1(\cB)\widetilde{m}^{c_2(\cB)}.
}
Substituting in our fixed value for $\varepsilon$, we obtain
\bes{
\bbE[X]\leq c_1(\cB)\widetilde{m}^{c_2(\cB)}m^{-\log(m)}+1
\lesssim c_1(\cB)m^{c_2(\cB)-\log(m)}+1
\leq C(\cB),
}
which proves the claim.
}
\prf{[Proof of Theorem~\ref{t:non-iid-sum}]
We follow the proof of Theorem \ref{t:non-iid-mix}, replacing Theorem \ref{t:main-prob-mixed} by Theorem \ref{t:main-prob-sum} and using $\varepsilon=1/2$. As there, we split into the cases $s\geq 2$ and $s=1$.

Assume first that $s\geq 2$. Then Theorem \ref{t:main-prob-sum} gives that, with probability at least $1/2$,
\be{
\label{eq:error_bound_for_det_thm_sum}
\nm{f-R((x_i,f(x_i))_{i=1}^m)}_{L^2}
\lesssim_{d,\beta,u}
\Big(\frac{1}{\widetilde{m}}\Big)^{g(\beta)}
\nm{f}_{H^{\beta}}
\quad
\forall f\in H^{\beta}(\bbT^d),\ g(\beta)>1/2.
}
This implies that there exists a choice of points $(x_1^{\star},\ldots,x_m^{\star})$ for which \eqref{eq:error_bound_for_det_thm_sum} holds with $x_i=x_i^{\star}$.

Now consider $s=1$. In this case we may select any points $(x_1^{\star},\ldots,x_m^{\star})\subset\bbT^d$. The argument in the proof of Theorem \ref{t:main-prob-sum} for $s=1$ does not use any condition on the sampling points and leads to the same bound \eqref{eq:error_bound_for_det_thm_sum}. This finishes the proof.}

\section{Lower bounds on nonlinear recovery widths}\label{s:lower-bounds-widths}

In this section, we present lower bounds for approximation in the spaces $H^{\alpha}_{\mathsf{mix}}$ and $H^{\beta}$. The main purpose of these results is to confirm that the algorithms established in Theorems \ref{t:main-prob-mixed}--\ref{t:non-iid-mix},\ref{t:main-prob-sum}--\ref{t:non-iid-sum} are optimal, up to the polylogarithmic factor arising in $\widetilde{m}$.
The following two theorems estimate the adaptive $m$-width \ef{adaptive-width}

\thm{\label{t:lower-bound-of-adaptive-m-widths-Halpha-mix} 
Let $\alpha=(\alpha_1,\dots,\alpha_d)$ with $\alpha > 1/2$ and $B(H^{\alpha}_{\mathsf{mix}})$ be the unit ball of $H^{\alpha}_{\mathsf{mix}}(\bbT^d)$. 
Then, for every $m\in\bbN$, the width $\varepsilon_m$ of $B(H^{\alpha}_{\mathsf{mix}})$ in $L^2$ satisfies
\begin{align*} 
\varepsilon_m(B(H^{\alpha}_{\mathsf{mix}}),L^2) \asymp_{\alpha,d} \left(\frac{(\log m)^{p(\alpha)-1}}{m}\right)^{h(\alpha)}. 
\end{align*}
}

\thm{\label{t:lower-bound-of-adaptive-m-widths-Hbeta} 
Let $\beta=(\beta_1,\dots,\beta_d)$ with $g(\beta)>1/2$ and $B(H^\beta)$ be the unit ball in $H^\beta(\bbT^d)$. 
Then, for every $m\in\bbN$, the width $\varepsilon_m$ of $B(H^{\beta})$ in $L^2$ satisfies
\begin{align*} 
\varepsilon_m(B(H^{\beta}),L^2) \asymp_{\beta,d}\left(\frac1m\right)^{g(\beta)}. 
\end{align*} 
}

Overall, these results confirm optimality of the algorithms developed in \S \ref{s:universal-algorithms-thms}, up to a polylogarithmic factor. 
Furthermore, by establishing a lower bound on $\varepsilon_m$, which allows for arbitrary adaptive linear measurements, they confirm that not only do pointwise samples constitute near-optimal information for recovery in $H^{\alpha}_{\mathsf{mix}}$ and $H^{\beta}$, but also that i.i.d.\ sampling from the underlying uniform measure is also near-optimal. In particular, adaptive sampling, or even a change of the probability measure employed for i.i.d.\ sampling, is unnecessary to achieve near-optimal rates. 
We remark in passing that the lower bound in Theorem \ref{t:lower-bound-of-adaptive-m-widths-Halpha-mix} was shown in \cite[Cor.\ 8.2]{byrenheid2017optimal} for general $L^p$ and $L^q$-spaces (our result corresponds to $p = q = 2$). However, the upper bound found in \cite[Thm.\ 8.4]{byrenheid2017optimal} only holds for $p \neq q$. Our proof is specific to $p = q = 2$ and follows different arguments.

As commented, the algorithms developed in \S \ref{s:universal-algorithms-thms} are universal, in that they achieve near-optimal rates simultaneously for different values of $\alpha$ and $\beta$. It is notable that the algorithms implied in these theorems only need work for a fixed value of $\alpha$ or $\beta$. It is an open problem whether a tighter lower bound can be established for universal algorithms, thereby narrowing the polylogarithmic gap between the current upper and lower bounds.

The remainder of this section establishes these results. We commence with the following, which demonstrates in an abstract sense that the adaptive $m$-width of an ellipse in $L^2(\bbT^d)$ is lower bounded by the worst-case best $s$-term approximation error.

\lem{\label{l:lower-bound-of-adaptive-m-widths}
Let $w=(w_n)_{n\in\bbZ^d} \in \ell^2(\bbZ^d)$ be a sequence of positive weights, and let \[ F = \left \{ f\in L^2(\bbT^d) : \sum_{n\in \bbZ^d}w_n^{-2}|\hat f_n|^2\le 1 \right \}. \]
Then, for every $m\in\mathbb{N}$, the width $\varepsilon_m(F,L^2)$ satisfies
\[ w^*_{2m+1}\le\varepsilon_m(F,L^2)  \le w^*_{m+1},
\] 
where $(w^*_i)^{\infty}_{i=1}$ is a non-increasing rearrangement of $w$.
}

\prf{
The proof follows along the lines of the proof of 
\cite[Theorem~5.4]{adcock2025sample}. 
We first require some notation.
For an index set $I$, sequence of positive weights $w=(w_i)_{i\in I}$ and $p$ such that $1\le p\leq\infty$, we write $\ell^p_w(I)$ for the set of real valued sequences $c=(c_i)_{i\in I}\in\bbR^{|I|}$ such that 
\[
\|c\|_{p,w} 
: = 
\begin{cases} \left(\sum_{n\in I} w_n^{-p}|c_n|^p\right)^{1/p} & 1 \leq p < \infty \\ \sup_{n\in I} w_n^{-1}|c_n| & p = \infty \end{cases} 
<\infty .
\] 
Let $B({\ell^p_w}(I))$ denote the unit ball in $\ell^p_w(I)$.  

Now let $S:C(\bbT^d)\to\mathbb{C}^m$ be an adaptive linear operator and 
$R:\mathbb{C}^m\to L^2(\bbT^d)$ be a reconstruction map. 
Define $J:\bbC^m\to\bbR^{2m}$ and $K:\bbR^{2m}\to\bbC^m$ by
$J(z)=(\Re z_1,\Im z_1,\dots,\Re z_m,\Im z_m)$ and
$K(a_1,b_1,\dots,a_m,b_m)=(a_1+\I b_1,\dots,a_m+\I b_m)$.
For a finite index set $I\subset\bbZ^d$, 
define the adaptive linear operator 
$M_I:\bbR^{|I|}\to\bbR^{2m}$ by
\[
M_I(c)= J\left(S\Bigl(\sum_{n\in I} c_n\phi_n\Bigr)\right),\quad \forall c = (c_n)_{n \in I} \in \bbR^{|I|}.
\]
Fix $c=(c_n)_{n\in I}\in B(\ell^2_w(I))$ and let $f=\sum_{n\in I}c_n\phi_n$.
Then
\begin{align*} 
\|f-R(S(f))\|_{L^2}^2 
& =\sum_{n\in\bbZ^d}\left|\frac{1}{(2\pi)^{d/2}}\int_{\bbT^d}\bigl(f-R(S(f))\bigr)(x)\phi_{-n}(x) d x\right|^2 \\ 
& \geq \sum_{n\in I}\left|c_n-\frac{1}{(2\pi)^{d/2}}\int_{\bbT^d}R(K(J(S(f))))(x)\phi_{-n}(x) d x\right|^2
\\
&=\sum_{n\in I}\left|c_n-\frac{1}{(2\pi)^{d/2}}\int_{\bbT^d}R(K(M_I(c)))(x)\phi_{-n}(x) d x\right|^2
\end{align*} 
Thus, defining  
$L:\bbR^{2m}\to\bbR^{|I|}$ as  
\[
L(z)=\Re\left(\frac{1}{(2\pi)^{d/2}}\int_{\bbT^d} R(K(z))(x)\phi_{-n}(x) \D x\right)_{n\in I},
\] 
we obtain
\begin{align*}
\|f-R(S(f))\|_{L^2}^2
\ge\sum_{n\in I}\left|c_n-\frac{1}{(2\pi)^{d/2}}\int_{\bbT^d}R(K(M_I(c)))(x)\phi_{-n}(x) d x\right|^2 \ge
\|c-L(M_I(c))\|_2^2.
\end{align*} 
We deduce that $\sup_{f\in \cF}\|f-R(S(f))\|_{L^2}^2\ge\sup_{c\in B({\ell^2_w}(I))}\|c-L(M_I(c))\|_{2}^2$.
Since $M_I$ is an adaptive linear operator, it follows that
\begin{align*}
\varepsilon_m(\cF,L^2)&=\inf\Bigl\{ 
\sup_{f\in \cF}\|f-R(S(f))\|_{L^2} 
:  
S:\cF\to\mathbb{C}^m \text{ adaptive}, 
R:\mathbb{C}^m\to L^2 
\Bigr\}
\\
& \ge \inf\Bigl\{ 
\sup_{c\in B({\ell^2_w}(I))}\|c-L(M_I(c))\|_{2} 
:  
M_I:\bbR^{|I|}\to\mathbb{R}^{2m} \text{ adaptive}, 
L:\mathbb{R}^{2m}\to \bbR^{|I|} 
\Bigr\}
\\
& = :
E_{2m}(B({\ell^2_w}(I)),\ell^2(I)). 
\end{align*} 
In the next step we invoke a standard result, which states that if $F$, subset of a normed space $(Z,\nm{\cdot}_Z)$, is symmetric with respect to the origin, i.e., $F=-F$, then $E_{m}(F,Z) 
\ge 
d^{m}(F,Z)$,
where
\bes{
d^{m}(F,Z) = \inf\Bigl\{ 
\sup_{z\in F\cap S}\|z\|_Z 
: 
S\subset Z \text{ with } \operatorname{codim}(S)\le m 
\Bigr\}
}
is the Gelfand width.
With this, and the fact that $B({\ell^2_w}(I))$ is symmetric, 
we have 
\begin{align}\label{eq:symmetric}
\varepsilon_{2m}(F,L^2) \ge E_{2m}(B({\ell^2_w}(I)),\ell^2(I)) 
\ge 
d^{2m}(B({\ell^2_w}(I)),\ell^2(I)). 
\end{align}
For the next step of our proof, we make use of a duality result. Consider the Kolmogorov width of a subset $F$ of a normed space $(Z,\nm{\cdot}_Z)$:
\[
d_m(F,Z) 
= 
\inf\Bigl\{ 
\sup_{f\in F}\inf_{z\in Z_m}\|f-z\|_Z 
: 
Z_m\subset Z \text{ with } \dim(Z_m)\le m 
\Bigr\}. 
\] 
Then
\cite[Theorem B.3]{adcock2024optimal} says that, given 
real numbers $p$ and $q$ such that $1\le p,q\le\infty$ and 
a finite sequence of positive weights $w=(w_n)_{n\in I}$, we have 
\begin{align*}
    d^m(B({\ell^{q^*}_{w}}(I)),\ell^{p^*}(I)) 
= 
d_m(B({\ell^p}(I)),\ell_{1/w}^q(I)). 
\end{align*}
where $p^*$ and $q^*$ satisfy $1/p^*+1/p=1$ and $1/q^*+1/q=1$ respectively. With this, we have
\begin{align}\label{eq:width-duality} 
d^{2m}(B({\ell^2_w}(I)),\ell^2(I)) 
= 
d_{2m}(B({\ell^{2}}(I)),\ell_{1/w}^2(I)). 
\end{align} 
We combine this with another technical result. 
\cite[Lemma B.4]{adcock2024optimal} states that given a finite sequence of positive weights $w=(w_n)_{n\in I}$ and $p,q$ such that $1\le p,q\le\infty$, we have 
\begin{align*}
    d_m(B({\ell^p}(I)),\ell^q_{1/w}(I)) 
= 
d_m(B(\ell^p_w(I)),\ell^q(I)). 
\end{align*}
With this, we have
\begin{align}\label{eq:width-reciprocal} 
    d_{2m}(B({\ell^2}(I)),\ell^2_{1/w}(I)) 
= 
d_{2m}(B({\ell^2_w}(I)),\ell^2(I)). 
\end{align} 
Combining \ef{eq:symmetric}, \ef{eq:width-duality} and \ef{eq:width-reciprocal} we deduce that $E_{2m}(F,L^2) 
\ge 
d_{2m}(B({\ell^2_w}(I)),\ell^2(I))$.  
This holds for any $I \subset \bbZ^d$.
Now, we make an explicit choice for $I$. Choose $N>2m$ and let $I=\{n_1,\dots,n_N\}$, where for $j\in\bbN$, $n_j\in \bbZ^d$ is the index of the $j$th largest entry of the sequence $w$. Fix $p>2$. 
By Hölder’s inequality, 
$
N^{-\frac{p-2}{2p}} B({\ell^p_w}(I))\subset B({\ell^2_w}(I)),
$ 
which implies that
\begin{align}\label{eq:holder-stesin}
d_{2m}(B({\ell^2_w}(I)),\ell^2(I)) 
\ge 
N^{-\frac{p-2}{2p}} d_{2m}(B({\ell^p_w}(I)),\ell^2(I)). 
\end{align}
In this next step, we invoke \cite[Theorem 3]{stesin1975aleksandrov}, which states that
for $N\in\bbN$ with $N>2m$, a finite sequence $w=(w_i)_{i\in I}$ of positive weights such that $|I|=N$ and $p,q$ such that $1\le q<p\le\infty$, 
we have
\[
d_m(B({\ell^p_w}(I)),\ell^q(I)) 
= 
\left( 
\min_{\substack{S\subset I\\|S|=N-m}} 
\sum_{n\in S} \left({w_n}\right)^{\frac{qp}{p-q}} 
\right)^{\frac{p-q}{qp}}. 
\] 
Combining this with \ef{eq:holder-stesin}, we have
\[
d_{2m}(B({\ell^2_w}(I)),\ell^2(I)) 
\ge 
N^{-\frac{p-2}{2p}}\left( 
\min_{\substack{S\subset I\\|S|=N-2m}} 
\sum_{n\in S}  \left({w_{n}}\right)^{\frac{2p}{p-2}} 
\right)^{\frac{p-2}{2p}} 
\ge
N^{-\frac{p-2}{2p}} w^{*}_{2m+1}. 
\] 
Letting $p\to2^+$ gives $\varepsilon_m(F,L^2) 
\ge w^*_{2m+1}$.

To complete the proof, it remains to show the upper bound. For $i = 1,2\ldots$, let $n_i \in \bbZ^d$ be such that $w^*_i = | w_{n_i} |$. Define the sampling-recovery pair $(S,R)$ by $S(f) =  ( \hat{f}_{n_i}  )^m_{i=1}$ and $R(c) = \sum^{m}_{i=1} c_i \phi_{n_i}$.
Then, for any $f \in F$,
\bes{
\nm{f - R(S(f))}^2_{L^2} = \sum_{i > m} | \hat{f}_{n_i}|^2 \leq (w^*_{m+1})^2 \sum_{i > m} (w^*_i)^{-2} | \hat{f}_{n_i}|^2 \leq (w^*_{m+1})^2.
}
We deduce that $\varepsilon_m(F,L^2)\leq w^*_{m+1}$,
as required.
}

\prf{[Proof of Theorem~\ref{t:lower-bound-of-adaptive-m-widths-Halpha-mix}] 
We invoke Lemma~\ref{l:lower-bound-of-adaptive-m-widths} with weights $w_n = \frac{1}{\prod_{j\in[d]} (1+|n_j|)^{\alpha_j}}$. Let $\varepsilon > 0$ and define 
$
S(\varepsilon) = \{ n \in \bbZ^d : w_n > \varepsilon \}.
$
First choose $\varepsilon_1$ such that $\varepsilon_1>0$ and $|S(\varepsilon_1)| \le m$. Then $w^*_{m+1} \le \varepsilon_1$ and Lemma~\ref{l:lower-bound-of-adaptive-m-widths} gives $\varepsilon_m(B(H^{\alpha}_{\mathsf{mix}}),L^2) \le w^*_{m+1} \le \varepsilon_1$.
Next choose $\varepsilon_2>0$ such that $|S(\varepsilon_2)| \geq 2m+1$. Then $w^*_{2m+1} \geq \varepsilon_2$ and Lemma~\ref{l:lower-bound-of-adaptive-m-widths} gives $\varepsilon_m(B(H^{\alpha}_{\mathsf{mix}}),L^2) \ge w^*_{2m+1} \ge \varepsilon_2$. 
By Lemma~\ref{l:adm} we have $|S(\varepsilon)| = |A(d,1/\varepsilon,\alpha)| \asymp_{d,\alpha} \varepsilon^{-\frac{1}{h(\alpha)}} (\log(1/\varepsilon))^{p(\alpha)-1}$. 
Therefore, there are constant $c_1 = c_1(d,\alpha),c_2 = c_2(d,\alpha)$ such that setting 
\begin{align*} 
\varepsilon_1 = c_{1}\left(\frac{\log(m)^{p(\alpha)-1}}{m}\right)^{h(\alpha)},\quad \varepsilon_2 = c_{2}\left(\frac{\log(2m)^{p(\alpha)-1}}{2m}\right)^{h(\alpha)} 
\end{align*} 
implies $|S(\varepsilon_1)|\le m$ and $|S(\varepsilon_2)|\ge 2m+1$. 
The result now follows.
}

\prf{[Proof of Theorem~\ref{t:lower-bound-of-adaptive-m-widths-Hbeta}]
We invoke Lemma~\ref{l:lower-bound-of-adaptive-m-widths} with weights given by $w_n = \frac{1}{1+\sum_{j\in[d]} |n_j|^{\beta_j}}$. We argue in the same way as in the previous proof, using Lemma \ref{lem:cdm} instead.
}

\section{Lower bounds for universal linear recovery}\label{s:universal-lower}
We now expose the gap between the performance of linear algorithms for universal recovery of anisotropic Sobolev functions and the nonlinear algorithms constructed in this work. We do this by following and modifying arguments found in \cite[\S 5.4]{temlyakov2018multivariate}, which are based on \cite{temlyakov1988approximation}. These works consider the spaces $H^{\beta}$, and therefore part of our effort is to extend them to the spaces $H^{\alpha}_{\mathsf{mix}}$. For simplicity, we only consider the case $p = q = 2$, while \cite{temlyakov1988approximation} considers arbitrary $1 \leq p,q\leq \infty$.

\subsection{Lower bounds for universal linear recovery}
We present the main results of this section. In the theorem below, we use the notation $\cL_m(B)$ to denote the set of linear operators $G$ whose domain contains all trigonometric polynomials, whose range is an $m$-dimensional subspace of $L^2(\bbT^d)$ and for which $\nm{G(\E^{\I n \cdot x })}_{L^2} \leq B$, $\forall n \in \bbZ^d$.

\thm{
[Lower bounds for universal linear recovery]
\label{t:main_theorem_index_of_universality}
Let $\cB \subset (0,\infty)^d$ have nonempty interior, and let $G \in \cL_m(B)$, for $m>m^*$ for some fixed natural number $m^*$ depending on $\cB$ and $d$ only, satisfy
\be{
\sup_{f \in B(H^{\beta})} \nm{f - G(f)}_{L^2} \leq C n^{-g(\beta)},\quad \forall \beta \in \cB,
}
where $C$ does not depend on $\beta$ or $n$. Then
\be{
m \gtrsim_{d,\cB,B,C} n (\log(n))^{d-1}.
}
Similarly, let $\cA \subset (0,\infty)^d$ have nonempty interior. If $G \in \cL_m(B)$, for $m>m^*$ for some fixed natural number $m^*$ depending on $\cA$ and $d$ only, satisfies
\be{
\sup_{f \in B(H^{\alpha}_{\mathsf{mix}})} \nm{f - G(f)}_{L^2} \leq C n^{-h(\alpha)} (\log(n))^{h(\alpha)(p(\alpha)-1)},\quad \forall \alpha \in \cA,
}
where $C$ does not depend on $\alpha$ or $n$, then
\be{
\label{pA-def}
m \gtrsim_{d,\cA,B,C} n (\log(n))^{p(\cA)-1},\qquad \text{where }
p(\cA) = \max_{\alpha\in \mathrm{int}(\cA)} p(\alpha).
}
}
The result above states that the loss of efficiency of the rate of universal linear recovery for either the $H^{\beta}$ or $H^{\alpha}_{\mathsf{mix}}$ spaces is always logarithmic in $n$. In particular, for the former, it behaves like $(\log(n))^{d-1}$ whenever the parameter set $\cB$ has nonempty interior. For the latter, it depends on $p(\cA)$. In particular, it also behaves like $(\log(n))^{d-1}$ whenever $\mathrm{int}(\cA)$ contains an element of the form $\alpha = (\alpha_0,\ldots,\alpha_0)$ for some $\alpha_0 > 0$.
This loss of efficiency of the universal recovery rate is connected to the concept of index of universality from \cite[\S 5.4]{temlyakov2018multivariate}.

\cor{
[Necessity of nonlinear algorithms for $H^{\beta}$]
\label{cor-necessity-Hbeta}
Let $\cB \subset (0,\infty)^d$ have nonempty interior and let $G \in \cL_m(B)$ for $m>m^*$ for some fixed natural number $m^*$ depending on $\cB$ and $d$ only,. Suppose that
\be{
\label{G-ub}
\sup_{f \in B(H^{\beta})} \nm{f - G(f)}_{L^2} \lesssim_{\cB,d} n^{-g(\beta)},\quad \forall \beta \in \cB.
}
Then
\bes{
n \lesssim_{\cB,B,d} \frac{m}{(\log m)^{d-1}}.
}
Conversely, let $H : C(\bbT^d) \to L^2(\bbT^d)$ be the algorithm from Theorem \ref{t:non-iid-sum}. Then
\bes{
\sup_{f \in B(H^{\beta})} \nm{f - H(f)}_{L^2} \lesssim_{\beta,d} n^{-g(\beta)},\quad \forall \beta \in \cB,
}
where $n \in \bbN$ satisfies
\bes{
n \ge \frac{m}{\log^3(m)\log(\log(m))}.
}
}

\cor{
[Necessity of nonlinear algorithms for $H^{\alpha}_{\mathsf{mix}}$]
\label{cor-necessity-Halphamix}
Let $\cA \subset (0,\infty)^d$ have nonempty interior and let $G \in \cL_m(B)$ for $m>m^*$ for some fixed natural number $m^*$ depending on $\cA$ and $d$ only. Suppose that
\be{
\sup_{f \in B(H^{\alpha}_{\mathsf{mix}})} \nm{f - G(f)}_{L^2}
\lesssim_{\cA,d} n^{-h(\alpha)} (\log(n))^{h(\alpha)(p(\alpha)-1)},\quad \forall \alpha \in \cA.
}
Then
\bes{
n \lesssim_{\cA,B,d} \frac{m}{(\log m)^{p(\cA)-1}},
}
where $p(\cA)$ is as in \ef{pA-def}. Conversely, let $H : C(\bbT^d) \to L^2(\bbT^d)$ be the algorithm from Theorem \ref{t:non-iid-mix}. Then
\bes{
\sup_{f \in B(H^{\alpha}_{\mathsf{mix}})} \nm{f - H(f)}_{L^2}
\lesssim_{\alpha,d} n^{-h(\alpha)} (\log(n))^{h(\alpha)(p(\alpha)-1)},\quad \forall \alpha \in \cA,
}
where $n \in \bbN$ satisfies
\bes{
n \ge \frac{m}{\log^3(m)\log(\log(m))}.
}
}

\prf{
[Proof of Corollaries \ref{cor-necessity-Hbeta} and \ref{cor-necessity-Halphamix}]
The first part follows directly from Theorem \ref{t:main_theorem_index_of_universality} by rearranging the bound on $m$. 
The second part follows from Theorem \ref{t:non-iid-sum} or Theorem \ref{t:non-iid-mix}.
}

Note that we could also use Theorems \ref{t:main-prob-mixed}, \ref{t:main-prob-sum} or \ref{t:main-expec-mix}, \ref{t:main-expec-sum} mutatis mutandis. Overall, these results confirm the necessity of nonlinear algorithms. In particular, linear algorithms incur a curse of dimensionality in the term $(\log m)^{d-1}$, while for nonlinear algorithms the corresponding term is independent of dimension and at most $\log^3(m) \log(\log(m))$. Note that for the $H^{\alpha}_{\mathsf{mix}}$ spaces, this holds whenever $p(\cA) = d$, i.e., when $\mathrm{int}(\cA)$ contains an element of the form $\alpha = (\alpha_0,\ldots,\alpha_0)$.

\subsection{Proof of Theorem \ref{t:main_theorem_index_of_universality}}

The proof follows from the following lemma, which is a modification of \cite[Thm. 5.4.3]{temlyakov2018multivariate}. 

\lem{\label{preliminary-result-for-Temlyakov-lower-bounds}
Let $0 < A_j < B_j$, $j = 1,\ldots,d$, $P = \prod^{d}_{j=1} [A_j , B_j ]$ and suppose that $G \in \cL_m(B)$, for $m>m^*$ for some fixed natural number $m^*$ depending on $P$ and $d$ only, satisfies either
\be{
\label{classical-temlyakov-cond}
\sup_{f \in B(H^{\beta})} \nm{f - G(f)}_{L^2} \leq C n^{-g(\beta)},\quad \forall \beta \in P,
}
or
\be{
\label{mixed-temlyakov-cond}
\sup_{f \in B(H^{\alpha}_{\mathsf{mix}})} \nm{f - G(f)}_{L^2} \leq C n^{-h(\alpha)} (\log(n))^{h(\alpha)(p(\alpha)-1)},\quad \forall \alpha \in P,
}
where $C$ does not depend on $\beta$ (respectively, $\alpha$) or $n$. In the latter case, suppose also that $\cap^{d}_{j=1} [A_j , B_j ) \neq \emptyset$.
Then
\bes{
m \gtrsim_{d,P,B,C} n (\log(n))^{d-1}.
}
}

\begin{foldme}
Note that \ef{classical-temlyakov-cond} is precisely the case studied in \cite[Thm.\ 5.4.3]{temlyakov2018multivariate}. That the result also holds under condition \ef{mixed-temlyakov-cond} is, to the best of our knowledge, new.
\end{foldme}

\prf{
We follow the construction in the proof of \cite[Thm.\ 5.4.3]{temlyakov2018multivariate}. This asserts the existence of a $\beta \in P$ (determined explicitly from $P$), natural numbers $k,s^*_1,\ldots,s^*_d$ and trigonometric polynomial $\nu$ with the following properties:
\begin{enumerate}[label=(\roman*)]
\item $2^k k^{d-1} \lesssim_{d,P,B} m$;
\item $\nm{\nu}_{L^1} = 1$;
\item $\hat{\nu}_n = 0$ if $|n_j| \geq 2^{s^*_j}$ for some $j$, where $s^*_j \in \bbN$ with $\sum^{d}_{j=1} s^*_j = k$;
\item $\nm{\nu - G(\nu)}_{L^2} \gtrsim_{d,P,B} 2^{k/2}$.
\item The vector $\beta^*$, defined explicitly by setting $ \beta^*_j = \frac{g(\beta) k}{s^*_j}$, $j = 1,\ldots,d$ belongs to $P$.
\end{enumerate}
Consider the anisotropic Sobolev spaces. By construction 
 $g(\beta^*) = g(\beta)$. Define the function
\bes{
f(x) = 2^{-k (g(\beta^*)+1/2)} \nu(x).
}
Then (iii) implies that $\nm{f}_{H^{\beta^{\star}}} \lesssim_{d,\beta^{*}} 2^{-k/2} \nm{\nu}_{L^2} \lesssim_{d,P} 2^{-k/2} \nm{\nu}_{L^2}$,
where the last inequality follows from the fact that $\beta$ is explicitly constructed from $P$. Further, Nikolskii's inequality (see, e.g., \cite[Thm. 3.3.2]{temlyakov2018multivariate}), (ii) and (iii) imply that
\bes{
\nm{f}_{H^{\beta^{\star}}} \lesssim_{d,P} 2^{-k/2} \prod^{d}_{j=1} (2^{s^{*}_j})^{1/2} \nm{\nu}_{L^1} = 1.
}
Hence $\nm{f}_{H^{\beta^*}} \lesssim_{d,P} 1$. Using (iv) and \ef{classical-temlyakov-cond}, we see that
\bes{
2^{-k g(\beta^*)} \lesssim_{d,P,B} \nm{f - G(f)}_{L^2} \lesssim_{d,P} C n^{-g(\beta^*)}
}
Thus $n \lesssim_{d,P,B,C} 2^k$, which implies that $n (\log(n))^{d-1} \lesssim_{d,P,B} 2^k k^{d-1}$. The result for the anisotropic Sobolev spaces now follows from (i).

Now consider the dominating mixed smoothness spaces. Since $\cap^{d}_{j=1} [A_j , B_j ) \neq \emptyset$ by assumption, there exists a point $\alpha_0$ such that $(\alpha_0,\alpha_0,\ldots,\alpha_0) \in P$ with $\alpha_0 < B_j$, $\forall j$. Now let
\bes{
\alpha_{\delta} = (\alpha_0 , \alpha_0 + \delta , \ldots , \alpha_0 + \delta ).
}
For all sufficiently small $\delta$, we have $\alpha_{\delta} \in P$. Now define
\bes{
f(x) = 2^{-s^*_1 \alpha_0 - \sum^{d}_{j=2} s^*_j (\alpha_0+\delta) - k/2} \nu(x).
}
Then (ii), (iii) and Nikolskii's inequality once more imply that
\bes{
\nm{f}_{H^{\alpha_{\delta}}_\mathsf{mix}} \lesssim_{d,P} 2^{-k/2} \nm{\nu}_{L^2} \lesssim 1.
}
We deduce from (iv) and \ef{mixed-temlyakov-cond} and the fact that $p(\alpha_{\delta}) = 1$, $h(\alpha_{\delta}) = \alpha_0$ that
\bes{
2^{-s^*_1 \alpha_0 - \sum^{d}_{j=1} s^*_j (\alpha_0+\delta)} \lesssim_{d,P,B} \nm{f - G(f) }_{L^2} \lesssim_{d,P} C n^{-\alpha_0}.
}
Letting $\delta \rightarrow 0^{+}$ and applying (iii) once more, we get $n \lesssim_{d,P,B,C} 2^k$. The result follows from (i).
}

\prf{[Proof of Theorem~\ref{t:main_theorem_index_of_universality}]
Since it has nonempty interior by assumption, $\cB$ must contain a rectangle $P = \prod^{d}_{j=1} [A_j , B_j]$ with $0 < A_j < B_j$, $\forall j$. 
Further, since $G$ satisfies \ef{G-ub} by assumption, we have that $G$ satisfies \ef{classical-temlyakov-cond} for the rectangle $P$ and with constant $C = C_{d,\cB}$.  
Consequently, we deduce from Lemma \ref{preliminary-result-for-Temlyakov-lower-bounds} that $m \gtrsim_{d,B,\cB} n (\log(n))^{d-1}$. The result now follows immediately.  

The argument for the anisotropic dominating mixed smoothness Sobolev spaces proceeds by following a dimension-reduction argument. Let $p = p(\cA)$ and choose $\alpha' \in \mathrm{int}(\cA)$ such that $p(\alpha') = p$.  
Let $J = \{j_1,\dots,j_p\} \subset \{1,\dots,d\}$ be the set of indices where $\alpha'$ attains its minimum, i.e., $\alpha'_{j_1}=\cdots=\alpha'_{j_p}=h(\alpha')=\min_{1\le j\le d}\alpha'_j$,
and define the parameter set  $\cA_J = \{\alpha_\mathsf{trunc}=(\alpha_{j_1},\dots,\alpha_{j_p}) : \alpha\in \cA\}\subset(0,\infty)^p$.  
Because $\alpha' \in \mathrm{int}(\cA)$, $\cA$ must contain a rectangle $P=\prod_{j=1}^d[A_j,B_j)$ with $0<A_j<B_j$ for all $j$ and $\cap_{i=1}^p[A_{j_i},B_{j_i})\neq \emptyset$ 
, the point $(\alpha'_{j_1},\dots,\alpha'_{j_p})=(h(\alpha'),\dots,h(\alpha'))$ must belong to $\mathrm{int}(\cA_J)$ and $\cA_J$ as well must contain a rectangle \[P_J= \{\alpha_\mathsf{trunc}=(\alpha_{j_1},\dots,\alpha_{j_p}) : \alpha\in P\}\subset(0,\infty)^p.\]  
We now work in dimension $p$ using three linear maps $H,Q,R$.  
Define $H:L^2(\bbT^d)\to L^2(\bbT^d)$ by  
\bes{  
Hf = \sum_{k\in\bbZ^d \,:\, k_j = 0, \forall j\notin J} \ip{f}{\phi_k}\phi_k.  
}
Define $Q:L^2(\bbT^p)\to L^2(\bbT^d)$ by $(Qf)(x_1,\dots,x_d) = f(x_{j_1},\dots,x_{j_p})$. 
Let $V_J \subset L^2(\bbT^d)$ be the span of $\{\phi_k : k_j = 0,\ \forall j\notin J\}$. 
Define $R:V_J\to L^2(\bbT^p)$ by  $(Rf)(x_1,\dots,x_p) = f(y_1,\dots,y_d)$,  
where $y_{j_i}=x_i$ for $i=1,\dots,p$ and $y_j\in\bbT$ is some arbitrary fixed value for $j\notin J$.
Since $f\in V_J$, the value of $f(y_1,\dots,y_d)$ does not depend on $y_j$ for $j\notin J$.  
Now define $G':L^2(\bbT^p)\to L^2(\bbT^p)$ by $G'  = R \circ H \circ G \circ Q$. We next observe that $\|Hf\|_{L^2(\bbT^d)}\le \|f\|_{L^2(\bbT^d)}$ for all $f\in L^2(\bbT^d)$, $QRf=f$ for all $f\in V_J$, 
\bes{  
\|Qf\|_{L^2(\bbT^d)}=(2\pi)^{\frac{d-p}{2}}\|f\|_{L^2(\bbT^p)},\quad f\in L^2(\bbT^p),  
}  
and
\bes{  
\|g\|_{L^2(\bbT^d)}=(2\pi)^{\frac{d-p}{2}}\|Rg\|_{L^2(\bbT^p)},\quad g\in V_J.  
}  
Finally, we notice that
\begin{align}\label{property-of-operator-Q}
\widehat{Qf}_{k} =
\begin{cases}
(2\pi)^{(d-p)/2}\widehat{f}_{k_J}, & k_j=0, \forall j\notin J,
\\0, & \text{otherwise},
\end{cases}\quad \forall k \in \bbZ^d,
\end{align}
where, for $k\in \bbZ^d$, $k_J$ represents the vector $(k_{j_1},\dots,k_{j_p})\in\bbZ^p$.

Now, since the range of $G$ has dimension at most $m$, the range of $G'$ also has dimension at most $m$.
Moreover, for $n\in\bbZ^p$ we have
\eas{ 
\|G'(\E^{\I n\cdot x})\|_{L^2(\bbT^p)} = (2\pi)^{\frac{p-d}{2}}\|H G Q(\E^{\I n\cdot x})\|_{L^2(\bbT^d)}&\le (2\pi)^{\frac{p-d}{2}}\|GQ(\E^{\I n\cdot x})\|_{L^2(\bbT^d)} 
\le (2\pi)^{\frac{p-d}{2}}B, 
}
where we used the fact that $\|Hf\|_{L^2}\le \|f\|_{L^2}$ for the second last inequality and the facts $G\in\cL_m(B)$ and $Q(\E^{\I n\cdot x}) = \E^{\I n\cdot x}$ for the last inequality. 
Thus, we conclude that $G'\in\cL_m(B)$ in dimension $p$. 

Next we bound $\sup_{f \in B(H^{\alpha_\mathsf{trunc}}_\mathsf{mix})(\bbT^p)} \nm{f - G'(f)}_{L^2(\bbT^p)}$ over all $\alpha_\mathsf{trunc}\in P_J$. 
Specifically, we show that there exists a constant $C$, independent of $\alpha_\mathsf{trunc}$ and $n$, such that  
\bes{  
\sup_{f \in B(H^{\alpha_\mathsf{trunc}}_\mathsf{mix})(\bbT^p)} \nm{f - G'(f)}_{L^2(\bbT^p)} \le C n^{-h(\alpha_\mathsf{trunc})} (\log(n))^{h(\alpha_\mathsf{trunc})(p(\alpha_\mathsf{trunc})-1)}, \quad \forall \alpha_\mathsf{trunc} \in P_J.
}
Pick any $\alpha\in P$ and consider $\alpha_\mathsf{trunc}\in P_J$. What follows will be true for all $\alpha\in P$ as $\alpha$ was picked arbitrarily. Let $f\in B(H^{\alpha_\mathsf{trunc}}_\mathsf{mix})(\bbT^p)$ and set $g=(2\pi)^{\frac{p-d}{2}}Qf$. 
By definition of $B(H^{\alpha_\mathsf{trunc}}_\mathsf{mix})(\bbT^p)$, we have
\bes{ 
\sum_{k \in \bbZ^p} \prod_{i=1}^p (1+|k_{j_i}|)^{2\alpha_{j_i}} |\hat f_k|^2 \le 1.
}
Applying \ef{property-of-operator-Q}, this is equivalent to
\bes{ 
\sum_{k \in \bbZ^d} \prod_{j=1}^d (1+|k_j|)^{2\alpha_j} |\hat g_k|^2 \le 1. 
}
Now, looking at the term $\nm{f - G'(f)}_{L^2(\bbT^p)},$ we obtain
\bes{
\nm{f - G'(f)}_{L^2(\bbT^p)}=\nm{f - RHGQ(f)}_{L^2(\bbT^p)}=(2\pi)^{\frac{p-d}{2}}\nm{Qf - QRHGQ(f)}_{L^2(\bbT^p)}=\nm{g - HG(g)}_{L^2(\bbT^d)}
}
where we used the identity $QRf=f$ for all $f\in V_J$ in the last equality.
Since $g \in V_J$, we have $g = H g$, and hence
\bes{ 
\nm{g - H G g}_{L^2(\bbT^d)} =\nm{Hg - H G g}_{L^2(\bbT^d)}\le \nm{g - G g}_{L^2(\bbT^d)}. 
}
Combining the above bounds, we obtain
\bes{ 
\sup_{f \in B(H^{\alpha_\mathsf{trunc}}_\mathsf{mix})(\bbT^p)} \nm{f - G'(f)}_{L^2(\bbT^p)} \le \sup_{g \in B(H^{\alpha}_{\mathsf{mix}}(\bbT^d))} \nm{g - G(g)}_{L^2(\bbT^d)}.}
The claimed bound now follows from the corresponding estimate for $G$ on $\bbT^d$:
\bes{
\sup_{f \in B(H^{\alpha_\mathsf{trunc}}_\mathsf{mix})(\bbT^p)} \nm{f - G'(f)}_{L^2(\bbT^p)} \lesssim_{d,\cA} n^{-h(\alpha)}(\log(n))^{h(\alpha)(p(\alpha)-1)},\quad \forall \alpha \in P.
}
For a small enough rectangle $P$, we have $h(\alpha_\mathsf{trunc})=h(\alpha)$ and $p(\alpha_\mathsf{trunc})=p(\alpha)$ for all $\alpha \in P$. Hence, we have
\bes{  
\sup_{f \in B(H^{\alpha_\mathsf{trunc}}_\mathsf{mix})(\bbT^p)} \nm{f - G'(f)}_{L^2(\bbT^p)} \lesssim_{d,\cA} n^{-h(\alpha_\mathsf{trunc})}(\log(n))^{h(\alpha_\mathsf{trunc})(p(\alpha_\mathsf{trunc})-1)},\quad \forall \alpha_\mathsf{trunc} \in P_J.  
}  
Hence, we deduce from Lemma \ref{preliminary-result-for-Temlyakov-lower-bounds} that $m \gtrsim_{d,B,\cA} n (\log(n))^{p-1}$.  
The result now follows.  
}

\section*{Acknowledgements}

BA acknowledges support from the Natural Sciences and Engineering Research
Council of Canada (NSERC) through grant RGPIN/2470-2021. BA \& AG acknowledge the support
of FRQ (Fonds de recherche du Qu\'ebec) – Nature et Technologies through grant 359708.

\bibliographystyle{abbrv} 
\small
\bibliography{TemplateRefs}

\end{document}